\documentclass[12pt]{article}
\usepackage{amssymb,amsmath}
\title{Studentized processes of $U$-statistics}
\author{Masoud M. Nasari\footnote{\footnotesize{Research supported by a Carleton university Faculty of Graduate Studies and Research scholarship, and NSERC Canada Discovery Grants of M. Cs\"{o}rg\H{o} and M. Mojirsheibani at Carleton university.}} \\
\footnotesize{\emph{School of Mathematics and Statistics, Carleton University,}}\\
\footnotesize{{\emph{ 1125 Colonel By Drive, Ottawa, Canada}}}\\
\footnotesize{\emph{e-mail}: \texttt{mmnasari@connect.carleton.ca}}}
\date{}
\begin{document}
\maketitle
\abstract{A uniform in probability approximation is established for Studentized processes of non degenerate $U$-statistics of order $m\geq 2$ in terms of a standard Wiener process. The classical condition that the second moment of kernel of the underlying U-statistic exists is relaxed to having $\frac{5}{3}$ moments. Furthermore, the conditional expectation of the kernel is only assumed to be in the domain of attraction of the normal law (instead of the classical two moment condition).}
\section{Introduction and Background}
Let $X_{1},X_{2},\ldots$, be a sequence of non-degenerate real-valued i.i.d. random variables with distribution $F$. Let
$h(X_{1},\ldots,X_{m})$, symmetric in its arguments, be a Borel-measurable real-valued kernel of order $m\geq 1$, and consider the
parameter $\theta=\int_{\mathbb{R}^m} h(x_{1},\ldots ,x_{m})\ dF(x_{1})\ldots dF(x_{m})<\infty.$ The corresponding $U$-statistic (cf. Serfling  \cite{serf}
or Hoeffding \cite{hof}) is\\
$$U_{n}={n \choose m}^{-1}
\sum_{C(n,m)}\ h(X_{i_{1}},\ldots,X_{i_{m}})=\ [n]^{-m}\sum_{C'(n,m)}\ h(X_{i_{1}},\ldots,X_{i_{m}}),$$
 where $m\leq n$,  $\sum_{C(n,m)}$ and $\sum_{C'(n,m)}$ respectively stand for summing over $C(n,m)=\{1\leq i_{1}<\ldots <i_{m}\leq n\}$ and $C\ '(n,m)=\{1\leq i_{1}\neq\ldots \neq i_{m}\leq n\}$ and $[n]^{-m}:=\displaystyle{\frac{(n-m)!}{n!}}$.
  For further  use throughout, we define
   $$\tilde{h}_{1}(x)=\mathbb{E}\textbf{(}h(X_{1},\ldots,X_{m})-\theta | X_{1} =x\textbf{)}.$$
\textbf{Definition}. A sequence $X,X_{1},X_{2},\ldots,$ of i.i.d. random variables is said to be in the domain of attraction of the normal law ($X\in$ $DAN$) if there exist sequences of constants $A_{n}$ and $B_{n}>0$ such that, as $n\rightarrow \infty,$
$$\frac{\sum_{i=1}^{n}X_{i}-A_{n}}{B_{n}}\longrightarrow_{d} N(0,1).$$
\textbf{Remark 1}. Furtherer to this definition of $DAN$, it is known that $A_{n}$ can be taken as $n\mathbb{E}(X)$ and $B_{n}=n^{1/2} \ell_{X}(n)$, where $\ell_{X}(n)$ is a slowly varying function at infinity (i.e., $\lim_{n\rightarrow \infty}
\frac{\ell_{X}(nk)}{\ell_{X}(n)}=1$ for any $k>0$), defined by the distribution of $X$. Moreover,
$\ell_{X}(n)=\sqrt{Var(X)}>0$, if $Var(X)< \infty$, and $\ell_{X}(n)\rightarrow \infty$, as $n\rightarrow \infty$, if $Var(X)=\infty$. Also $X$ has all moments less than 2, and the variance of $X$ is positive, but need not be finite.
\par
Noting that $\tilde{h}_{1}(X_{1}),\tilde{h}_{1}(X_{2}),\ldots,$ are i.i.d. random variables with mean zero $(\mathbb{E}\tilde{h}_{1}(X_{1})=0)$, Nasari (cf. \cite{nasa} ) observed that Proposition 2.1 of Cs\"{o}rg\H{o}, Szyszkowicz and Wang [CsSzW] (2004 \cite{cs4}) (cf. also Theorem 1 of [CsSzW] 2003 \cite{cs3}) reads as follows (cf. Lemma 2 in Nasari \cite{nasa})
\\
\textbf{Lemma A }. \emph{As} $n\rightarrow \infty$,  \emph{the following statements are equivalent}: \\ \\
(a)\ \ $\tilde{h}_{1}(X_{1})\in$ $DAN$; \\

\emph{There is a sequence of constants} $B_{n}\nearrow \infty$, \emph{such that}\\ \\
(b) $\ \displaystyle{\frac{\sum_{i=1}^{[nt_{0}]}\tilde{h}_{1}(X_{i})}{B_{n}}}\longrightarrow_{d} N(0,t_{0})$ \emph{for} $t_{0}\in(0,1];$ \\\\ \\
$(\textrm{c}) \ \displaystyle{\frac{\sum_{i=1}^{[nt]}\tilde{h}_{1}(X_{i})}{B_{n}}}\longrightarrow_{d} W(t)$ \emph{on} $(D[0,1],\rho)$,
\emph{where} $\rho$ \emph{is the sup-norm metric}

\emph{for functions in} $D[0,1]$ \emph{and} $\{W(t),0\leq t\leq 1\}$
\emph{is a standard Wiener}

\emph{process};\\\\
(d) \emph{On an appropriate probability space for} $X_{1},X_{2},\ldots,$ \emph{we can construct a}

\emph{standard Wiener process} $\{W(t),0\leq t<\infty\}$ \emph{such that}\\
$$\sup_{0\leq t \leq 1}\left|\frac{\sum_{i=1}^{[nt]}\tilde{h}_{1}(X_{i})}{B_{n}} -\ \frac{W(nt)}{n^{\frac{1}{2}}} \right|=o_{P}(1).$$
\\
\par
Here and throughout, $B_{n}$ is as in Remark 1, from now on written as $B_{n}=n^{1/2}\ell(n)$, where $\ell(.)$, the slowly varying function at infinity, is defined by the distribution of the random variable $\tilde{h}_{1}(X_{1})$ (cf. Remark 1).
\\
\textbf{Remark 2}. The statement (c), whose notion will be used throughout, stands for the following functional central limit theorem (cf. Remark 2.1 in Cs\"{o}rg\H{o}, Szyszkowicz and Wang [CsSzw] (2004) \cite{cs4}). On account of (d), as $n\rightarrow\infty$, we have
$$g(S_{[n.]}/V_{n})\longrightarrow_{d} g(W(.))$$
for all $g:D=D[0,1]\longrightarrow \mathbb{R}$ that are $(D,\mathfrak{D})$ measurable and  $\rho$-continuous, or $\rho$-continuous except at points forming a set of Wiener measure zero on $(D,\mathfrak{D})$, where $\mathfrak{D}$ denotes the $\sigma$-field of subsets of $D$ generated by the finite-dimensional subsets of $D$.
\\
\par
In view of (b) of Lemma A with $t_{0}=1$, Corollary 2.1 of [CsSzW] (2004 \cite{cs4}), i.e., Raikov's theorem as stated and proved in Gin\'{e}, G\"{o}tze and Mason (1997 \cite{gin2}), yields the following version of it in the present context.
\\
\textbf{Corollary A}. \emph{As} $n\rightarrow\infty$, \emph{we have}
$$\frac{1}{n\ell^{2}(n)}\ \sum_{i=1}^{n} \tilde{h}^{2}_{1}(X_{i})\longrightarrow_{P} 1.$$
\par
Nasari \cite{nasa} proved a projection approximation of $U_{n}$into sums of the i.i.d. random variables $\tilde{h}_{1}(X_{1}),\tilde{h}_{1}(X_{2}),\ldots,$ that reads as follows (cf. Theorem 3 of \cite{nasa}).
\\ \\
\textbf{Theorem A}.  \emph{If} $\mathbb{E}\ [ |h(X_{1},\dots,X_{m})|^{\frac{4}{3}}\log|h(X_{1},\ldots,X_{m})|\ ]<\infty$ \emph{and}  $\tilde{h}_{1} (X_{1})\in$

$DAN$, \emph{then, as} $n\rightarrow \infty$, \emph{we have}
\\ \\
  $$\sup_{0\leq t \leq 1} \left|\displaystyle{\frac{[nt]}{m}} \frac{U_{[nt]}-\theta}{B_{n}}- \displaystyle{\frac{\sum_{i=1}^{[nt]}\tilde{h}_{1}(X_{i})}{B_{n}}} \right|=o_{P}(1).$$
\\
\par
In view of Lemma A and Theorem A , Nasari \cite{nasa} concluded his Theorem 2 that reads as follows.
\\ \\
\textbf{Theorem B.} \emph{If}
\\ \\
(a) $\mathbb{E}\ [ |h(X_{1},\dots,X_{m})|^{\frac{4}{3}}\log|h(X_{1},\ldots,X_{m})|\ ]<\infty$  \emph{and}
$\tilde{h}_{1} (X_{1}) \in$ $DAN$,\\

\emph{then}, \emph{as} $n\rightarrow \infty$, \emph{we have}\\ \\
(b) $\displaystyle{\frac{[nt_{0}]}{m}\ \frac{U_{[nt_{0}]}-\theta }{B_{n}}}\longrightarrow_{d} \ N(0,t_{0}),\ where \ t_{0}\in(0,1]$;\\ \\
(c) $\displaystyle{\frac{[nt]}{m}\ \frac{U_{[nt]}-\theta }{B_{n}}}\ \longrightarrow_{d}\ W(t)$ \emph{on} ($D$[0,1],$\rho$), \emph{where} $\rho$ \emph{is the sup-norm for} \\

\emph{functions in}   $D$[0,1] \emph{and} $\{W(t), 0\leq t \leq 1\}$ \emph{is a standard Wiener process};\\ \\
(d) \emph{On an appropriate probability space for} $X_{1},X_{2},\ldots$, \emph{we can construct a}

\emph{standard Wiener process} $\{W(t), 0\leq t < \infty\}$ \emph{such that}
$$\sup_{0\leq t \leq 1}\ \left|\ \frac{[nt]}{m}\frac{U_{[nt]}-\theta }{B_{n}} -\ \frac{W(nt)}{n^{\frac{1}{2}}} \right|=o_{P}(1).$$
\\
\par
We note in passing that the weak convergence result of part (c) of Theorem B for non degenerate $U$-statistics extend those obtained by Miller and Sen in 1972 (cf. Theorem 1 of \cite{mil} )

\par
Define  the pseudo-selfnormalized $U$-process $U^{*}_{[nt]}$ as follows
\\ $$U_{[nt]}^{*}=\left\{
\begin{array}{ll}
\ 0 \qquad \qquad  \  ,& \hbox{$0\leq t <\displaystyle{\frac{m}{n}}, $}\\
                     \displaystyle{\frac{U_{[nt]}-\theta }{V_{n}}} \  \ \ \  , & \hbox{$\displaystyle{\frac{m}{n}}\leq                      t\leq 1,$}\\
                    \end{array}
                  \right.$$\\ where [.] denotes the greatest integer function and  $V_{n}^{2}:=\sum_{i=1}^{n}\tilde{h}^{2}_{1}(X_{i})$. Combining Theorem A with Corollary A,    Nasari (cf. \cite{nasa}) inferred his Theorem 1 which reads as follows.
\\ \\
\textbf{Theorem C.} $If$ \\ \\
(a) $\mathbb{E}[\ |h(X_{1},\dots,X_{m})|^{\frac{4}{3}}\log|h(X_{1},\ldots,X_{m})|\ ]<\infty \ and$
$\tilde{h}_{1} (X_{1}) \in DAN$,\\

$then,\ as\ n\rightarrow \infty,\ we \ have$\\ \\
(b) $\displaystyle{\frac{[nt_{0}]}{m}\ U^{*}_{[nt_{0}]}}\rightarrow_{d} \ N(0,t_{0}), \ for \ t_{0}\in(0,1]$;\\ \\
(c) $\displaystyle{\frac{[nt]}{m}\ U^{*}_{[nt]}}\ \rightarrow_{d}\ W(t)$ \emph{on} ($D$[0,1],$\rho$), $where$ $\rho$ \emph{is the sup-norm for functions in}\\

 $D[0,1] \ and\  \{W(t), 0\leq t \leq 1\}\ is\ a\ standard\ Wiener\ process$;\\ \\
(d) $On\ an\ appropriate\ probability \ space \ for\ X_{1},X_{2},\ldots,\ we\ can\ construct\ a$

$standard \ Wiener \ process\ \{W(t), 0\leq t < \infty\}\ such\ that$
$$\sup_{0\leq t \leq 1}\ \left|\ \frac{[nt]}{m}\ U^{*}_{[nt]} -\ \frac{W(nt)}{n^{\frac{1}{2}}} \right|=o_{P}(1).$$
\par
We note that in the light of Corollary A, a similarly pseudo-selfnormalized version of Lemma A is also immediate (cf. Lemma 1 in Nasari \cite{nasa} ). Moreover, these two lemmas, i.e., Lemmas 1 and 2 in Nasari \cite{nasa}, respectively  coincide with  Theorem 1 of [CsSzW] (2003 \cite{cs3}), and  with Proposition 2.1 of [CsSzW] (2004 \cite{cs4}). Thus Theorems B and C with $m\geq2$ amount to begin extensions of Theorem 1 of [CsSzW] (2003 \cite{cs3}) to $U$-statistics of order $m\geq2$.
\par
While, in view of Raikov's theorem as in Corollary A, Theorems B and C are equivalent, Theorem C as stated constitutes a significant first step toward studentizing $U$-statistics for the sake of establishing asymptotic confidence intervals for $\theta$ in a nonparametric manner (cf. Theorem 1 and Corollary 1 of the next session that, in turn, leads to Main Theorem of this exposition). The pseudo-selfnormalizing sequence $V_{n}$ of Theorem C still depends on the distribution function $F$ that can not usually assumed to be known. Hence our Theorem 1 in this exposition.

\section{Statement of the results}
For $i=1,\ldots,n$, let $U^{i}_{n-1}$ be the \emph{jackknifed} version of $U_{n}$ based on $X_{1},\ldots,X_{i-1},$ \\
$X_{i+1},\ldots,X_{n}$, defined  as follows.
$$U^{i}_{n-1}=\displaystyle{\frac{1}{{n-1 \choose m}}}
\sum_{\substack{1\leq j_{1}<\ldots<j_{m}\leq n \\ j_{1},\ldots, j_{m}\neq i}} h(X_{j_{1}},\ldots,X_{j_{m}}).$$
Also define the \emph{Studentized} $U$-process as follows.\\
$$\displaystyle{U}_{[nt]}^{\emph{stu}}=\left\{
\begin{array}{ll}
\ 0,& \hbox{$0\leq t <\displaystyle{\frac{m}{n}}, $}\\
                     \displaystyle{\frac{U_{[nt]}-\theta }{\sqrt{(n-1)\sum_{i=1}^{n}(U^{i}_{n-1}-U_{n})^{2}}}}, & \hbox{$\displaystyle{\frac{m}{n}}\leq t\leq 1.$}
                     \end{array}
                     \right.$$\\
\textbf{Remark 3.} Unlike the $U$-processes in Theorems B and C, apart from the parameter $\theta$ of interest, $\displaystyle{U}_{[nt]}^{\emph{stu}}$ is completely computable, based on the observations $X_{1},\ldots,X_{n}$.\\
\par
Under a slightly stronger moment condition, which is the price we pay for the normalization involved in $\displaystyle{U}_{[nt]}^{\emph{stu}}$,  the  Studentized companion of Theorems B and C reads  as follows.
\\ \\
\textbf{Main Theorem}. \emph{If}
\\ \\
(a) $\mathbb{E}|h(X_{1},\dots,X_{m})|^{\frac{5}{3}}<\infty$  \emph{and}
$\tilde{h}_{1} (X_{1}) \in$ $DAN$,\\

 $then,\ as\ n\rightarrow \infty,\ we \ have$\\ \\
(b) $\displaystyle{[nt_{0}]\ U^{\emph{stu}}_{[nt_{0}]}}\rightarrow_{d} \ N(0,t_{0}), \ for \ t_{0}\in(0,1]$;\\ \\
(c) $\displaystyle{[nt]\ U^{\emph{stu}}_{[nt]}} \rightarrow_{d}\ W(t)$ \emph{on} ($D$[0,1],$\rho$), $where$ $\rho$ \emph{is the sup-norm for functions in}\\

 $D[0,1] \ and\  \{W(t), 0\leq t \leq 1\}\ is\ a\ standard\ Wiener\ process$;\\ \\
(d) $On\ an\ appropriate\ probability \ space \ for\ X_{1},X_{2},\ldots,\ we\ can\ construct\ a$

$standard \ Wiener \ process\ \{W(t), 0\leq t < \infty\}\ such\ that$
$$\sup_{0\leq t \leq 1}\ \left|\ \displaystyle{[nt]\ U^{\emph{stu}}_{[nt]}} -\ \frac{W(nt)}{n^{\frac{1}{2}}} \right|=o_{P}(1).$$\\
\par
In view of Theorems B and C and on account of Raikov's theorem (cf. Corollary A), which via (b) of Lemma A with $t_{0}=1$ in this context states that, as $n\rightarrow\infty$,   $\displaystyle{\frac{1}{n\ \ell^{2}(n)} \sum_{i=1}^{n} \tilde{h}^{2}_{1}(X_{i})}\rightarrow_{P}1$, in order to prove Main Theorem it suffices to prove the following result.
\\\\
\textbf{Theorem 1}. \emph{If}  $\mathbb{E}|h(X_{1},\dots,X_{m})|^{\frac{5}{3}}<\infty$ $and$ $\tilde{h}_{1} (X_{1}) \in$ $DAN$, \emph{then}, \emph{as} $n\rightarrow\infty$,\\ \\
$$\left|\displaystyle{\frac{(n-1)}{m^{2}\ \ell^{2}(n)}} \ \sum_{i=1}^{n} (U^{i}_{n-1}-U_{n})^{2}-\displaystyle{\frac{1}{n\ \ell^{2}(n)} \sum_{i=1}^{n} \tilde{h}^{2}_{1}(X_{i})}\right|=o_{P}(1).$$\\
\par
 Consequently, the latter approximation combined with Corollary A yields  a Raikov type result for the distribution free jackkifed version of $U$-statistics which is of interest on its own (cf. Remark 4).
 \\ \\
\textbf{Corollary 1}. \emph{If}  $\mathbb{E}|h(X_{1},\dots,X_{m})|^{\frac{5}{3}}<\infty$ $and$ $\tilde{h}_{1} (X_{1}) \in$ $DAN$, \emph{then}, \emph{as} $n\rightarrow\infty$,\\ \\
$$\displaystyle{\frac{(n-1)}{m^{2}\ \ell^{2}(n)}} \ \sum_{i=1}^{n} (U^{i}_{n-1}-U_{n})^{2}\longrightarrow_{P}1.$$
\\ \\

Combining now Corollary 1 with Theorem B we arrive at Main Theorem of this paper.
\\
\textbf{Remark 4}. When $\mathbb{E}\ h^{2}(X_{1},\dots,X_{m})<\infty$, which in turn implies that $\mathbb{E}\tilde{h}_{1}^{2}(X_{1})<\infty$, then $\ell^{2}(n)=\mathbb{E}\tilde{h}_{1}^{2}(X_{1})>0$  and, as $n\rightarrow\infty$, Corollary 1 implies that\\
$$\displaystyle{\frac{(n-1)}{m^{2}\ }} \ \sum_{i=1}^{n} (U^{i}_{n-1}-U_{n})^{2}\longrightarrow_{P}\mathbb{E}\tilde{h}_{1}^{2}(X_{1}).$$
The latter version of Corollary 1 coincides with one of the result obtained by Arvesen \cite{Arv} who extended the idea of the so-called (by Tukey) pseudo- values to $U$-statistics and studied the asymptotic distribution of non-degenerate $U$-statistics via jackknifing.
\\
\textbf{Remark 5}. When $m=1$, the projection $\tilde{h}_{1}(X_{1})$ will coincide with $h(X_{1})-\theta$, then Main Theorem corresponds to Corollary 5 of [CsSzW] (2008 \cite{cs8a}) on taking the weight function $q=1$ for the therein studied Studentized process $T_{n,t}(X-\mu)$, i.e., when $m=1$, then the studentized $U$-process $U^{\emph{stu}}_{[nt]}$ coincides  with $T_{n,t}(X-\mu)$.  Hence in this exposition we shall state our proofs for $m\geq2.$  Also when $m=2$, the two conditions in (a) of Main Theorem as well as the idea of its proof by truncation, coincide with the corresponding ones of Theorem 2 of [CsSzW] (2008b \cite{cs8b}) on weighted approximations for Studentized $U$-type processes .

\section{Proofs}
To prove Theorem 1, it suffices to show that as $n\rightarrow\infty$,
$$\left|(n-1) \ \sum_{i=1}^{n} (U^{i}_{n-1}-U_{n})^{2}-\displaystyle{\frac{m^{2}}{n}} \sum_{i=1}^{n} \tilde{h}^{2}_{1}(X_{i})\right|=o_{P}(1).\qquad\qquad\qquad \qquad \qquad \qquad (1)$$\\\\
Before proving (1) we do some simplifications as follows.\\ \\
$(n-1) \ \sum_{i=1}^{n} (U^{i}_{n-1}-U_{n})^{2}$\\
$$=(n-1)\sum_{i=1}^{n}\left[\displaystyle{\frac{\binom{n}{m}}{\binom{n-1}{m}}} \ U_{n}-\ \displaystyle{\frac{1} {\binom{n-1}{m}}}\sum_{\substack{1\leq j_{1}<\ldots<j_{m-1}\leq n \\ j_{1},\ldots, j_{m-1
}\neq i}} h(X_{i},X_{j_{1}}\ldots,X_{j_{m-1}})-\ U_{n}\right]^{2}\ \ \ $$
$$=(n-1)\sum_{i=1}^{n}\left[\displaystyle{\frac{1} {\binom{n-1}{m}}}\sum_{\substack{1\leq j_{1}<\ldots<j_{m-1}\leq n \\ j_{1},\ldots, j_{m-1}\neq i}} h(X_{i},X_{j_{1}}\ldots,X_{j_{m-1}})-\ \textbf{(}\displaystyle{\frac{\binom{n}{m}}{\binom{n-1}{m}}}-1\textbf{)}\ U_{n}\right]^{2}  \qquad $$
$$=(n-1) \sum_{i_{1}=1}^{n} \left[ \displaystyle{\frac{m}{n-m}} \textbf{(}\displaystyle{\frac{1}{\binom{n-1}{m-1}}}\sum_{\substack{1\leq i_{2}<\ldots< i_{m}\leq n\\ i_{2},\ldots, i_{m}\neq i_{1}}} h(X_{i_{1}},X_{i_{2}},\ldots,X_{i_{m}})-\ U_{n}\textbf{)} \right]^{2}\qquad \qquad$$
$$=\displaystyle{\frac{m^{2}(n-1)}{(n-m)^{2}}} \sum_{i_{1}=1}^{n} \left[ \displaystyle{\frac{1}{\binom{n-1}{m-1}}}\sum_{\substack{1\leq i_{2}<\ldots< i_{m}\leq n\\ i_{2},\ldots, i_{m}\neq i_{1}}} h(X_{i_{1}},X_{i_{2}},\ldots,X_{i_{m}})-\ U_{n} \right]^{2}\qquad \qquad\ \ \ (\large{*})$$\\
$$=\displaystyle{\frac{m^{2}(n-1)}{(n-m)^{2}}} \sum_{i_{1}=1}^{n} \left[ \displaystyle{\frac{1}{\binom{n-1}{m-1}}}\sum_{\substack{1\leq i_{2}<\ldots< i_{m}\leq n\\ i_{2},\ldots, i_{m}\neq i_{1}}} h(X_{i_{1}},X_{i_{2}},\ldots,X_{i_{m}})\right]^{2}+\displaystyle{\frac{m^{2}n (n-1)}{(n-m)^{2}}}\ U_{n}^{2}$$\\
$$-\ 2 \ \displaystyle{\frac{m^{2}(n-1)}{(n-m)^{2}}}\ U_{n} \ \displaystyle{\frac{1}{\binom{n-1}{m-1}}} \ \sum_{i_{1}=1}^{n} \ \sum_{\substack{1\leq i_{2}<\ldots< i_{m}\leq n\\ i_{2},\ldots, i_{m}\neq i_{1}}}   h(X_{i_{1}},X_{i_{2}},\ldots,X_{i_{m}})\qquad \qquad \quad  $$\\
$$=\displaystyle{\frac{m^{2}(n-1)}{(n-m)^{2}}} \sum_{i_{1}=1}^{n} \left[ \displaystyle{\frac{1}{\binom{n-1}{m-1}}}\sum_{\substack{1\leq i_{2}<\ldots< i_{m}\leq n\\ i_{2},\ldots, i_{m}\neq i_{1}}} h(X_{i_{1}},X_{i_{2}},\ldots,X_{i_{m}})\right]^{2}+\displaystyle{\frac{m^{2}n (n-1)}{(n-m)^{2}}}\ U_{n}^{2}$$\\
$$-\ 2 \ \displaystyle{\frac{m^{2}(n-1)}{(n-m)^{2}}}\ U_{n} \ \displaystyle{\frac{1}{(m-1)!\ \binom{n-1}{m-1}}} \ \sum_{i_{1}=1}^{n} \ \sum_{\substack{1\leq i_{2}\ \neq\ldots\neq \ i_{m}\leq n\\ i_{2},\ldots, i_{m}\neq i_{1}}}   h(X_{i_{1}},X_{i_{2}},\ldots,X_{i_{m}})\qquad $$\\
$$=\displaystyle{\frac{m^{2}(n-1)}{(n-m)^{2}}} \sum_{i_{1}=1}^{n} \left[ \displaystyle{\frac{1}{\binom{n-1}{m-1}}}\sum_{\substack{1\leq i_{2}<\ldots< i_{m}\leq n\\ i_{2},\ldots, i_{m}\neq i_{1}}} h(X_{i_{1}},X_{i_{2}},\ldots,X_{i_{m}})\right]^{2}+\displaystyle{\frac{m^{2}n (n-1)}{(n-m)^{2}}}\ U_{n}^{2}$$\\
$$-\ 2 \ \displaystyle{\frac{m^{2}(n-1)}{(n-m)^{2}}}\ U_{n} \ \displaystyle{\frac{1}{(m-1)!\ \binom{n-1}{m-1}}} \ \sum_{1\leq i_{1}\neq i_{2}\ \neq\ldots\neq \ i_{m}\leq n }   h(X_{i_{1}},X_{i_{2}},\ldots,X_{i_{m}})\qquad$$\\
$$=\displaystyle{\frac{m^{2}(n-1)}{(n-m)^{2}}} \sum_{i_{1}=1}^{n} \left[ \displaystyle{\frac{1}{\binom{n-1}{m-1}}}\sum_{\substack{1\leq i_{2}<\ldots< i_{m}\leq n\\ i_{2},\ldots, i_{m}\neq i_{1}}} h(X_{i_{1}},X_{i_{2}},\ldots,X_{i_{m}})\right]^{2}+\displaystyle{\frac{m^{2}n (n-1)}{(n-m)^{2}}}\ U_{n}^{2}$$\\
-$\ 2 \ \displaystyle{\frac{m^{2}n (n-1)}{(n-m)^{2}}}\ U_{n}^{2}$\\
$$=\ \displaystyle{\frac{m^{2}(n-1)}{(n-m)^{2}}} \sum_{i_{1}=1}^{n} \left[ \displaystyle{\frac{1}{\binom{n-1}{m-1}}}\sum_{\substack{1\leq i_{2}<\ldots< i_{m}\leq n\\ i_{2},\ldots, i_{m}\neq i_{1}}} h(X_{i_{1}},X_{i_{2}},\ldots,X_{i_{m}})\right]^{2}\ -\ \displaystyle{\frac{m^{2}n (n-1)}{(n-m)^{2}}}\ U_{n}^{2} $$
$$\qquad \qquad \qquad \qquad \qquad \qquad \qquad \qquad \qquad \qquad \qquad \qquad \qquad \qquad \qquad \qquad \qquad (2)$$\\
\textbf{Remark 6}. In  view of (*)  in what will follow  without loss of generality we may and shall assume that $\theta=0$.\\
\par
In view of (2) to prove (1) it will be enough to prove the following two propositions.\\\\
\textbf{Proposition 1}. \emph{If} $\mathbb{E}|h(X_{1},\ldots,X_{m})|^{\frac{5}{3}}<\infty$, \emph{then, as} $n\rightarrow\infty$,
$$U_{n}^{2}\longrightarrow 0\ \ a.s.$$\\
\textit{\textbf{Proof of Proposition 1}}\\\\
The proof this theorem follows from the SLLN for $U$-statistics (cf. for example Serfling \cite{serf}).\\\\
\textbf{Proposition 2}. \emph{If} $\mathbb{E}|h(X_{1},\ldots,X_{m})|^{\frac{5}{3}}<\infty$ and $\tilde{h}_{1}(X_{1})\in DAN$ \emph{then, as} $n\rightarrow\infty$,\\ \\
$$\left|\frac{(n-1)}{(n-m)^{2}}\sum_{i_{1}=1}^{n} \left[ \displaystyle{\frac{1}{\binom{n-1}{m-1}}}\sum_{\substack{1\leq i_{2}<\ldots< i_{m}\leq n\\ i_{2},\ldots, i_{m}\neq i_{1}}} h(X_{i_{1}},X_{i_{2}},\ldots,X_{i_{m}})\right]^{2}-\frac{1}{n}\sum_{i=1}^{n}\tilde{h}^{2}_{1}(X_{i})\right|$$
\\
$$\qquad \qquad \qquad \qquad \qquad \qquad \qquad \qquad \qquad \qquad \qquad \qquad \qquad \qquad =o_{P}(1).$$\\\\
\textit{\textbf{Proof of Proposition 2}}\\
In what will follow  $a_{n}\thickapprox b_{n}$ stands for the asymptotic equivalency of the numerical sequences $(a_{n})_{n}$ and $(b_{n})_{n}$, i.e., as $n\rightarrow \infty$, $\displaystyle{\frac{a_{n}}{b_{n}}\rightarrow 1}$.
\par
To prove Proposition 2 observe that\\\\
$\displaystyle{\frac{(n-1)}{(n-m)^{2}}}\sum_{i_{1}=1}^{n} \left[ \displaystyle{\frac{1}{\binom{n-1}{m-1}}}\sum_{\substack{1\leq i_{2}<\ldots< i_{m}\leq n\\ i_{2},\ldots, i_{m}\neq i_{1}}} h(X_{i_{1}},X_{i_{2}},\ldots,X_{i_{m}})\right]^{2}$\\\\
$=  \displaystyle{\frac{(n-1)}{(n-m)^{2}}} \sum_{i_{1}=1}^{n} \left[ [n-1]^{-m+1}\sum_{\substack{1\leq i_{2}\neq\ldots\neq i_{m}\leq n\\ i_{2},\ldots, i_{m}\neq i_{1}}} h(X_{i_{1}},X_{i_{2}},\ldots,X_{i_{m}})\right]^{2}$\\
$$\thickapprox [n]^{-2m+1} \ \sum_{i_{1}=1}^{n} \left[ \sum_{\substack{1\leq i_{2}\neq\ldots \neq i_{m}\leq n\\ i_{2},\ldots, i_{m}\neq i_{1}}} h(X_{i_{1}},X_{i_{2}},\ldots,X_{i_{m}}) \right]^{2}\qquad \qquad \qquad \qquad \qquad \quad \quad $$\\
$$\textbf{=}\  [n]^{-2m+1} \sum_{1\leq i_{1}\neq\ldots\neq i_{m}\leq n}h^{2}(X_{i_{1}},\ldots,X_{i_{m}}) \qquad \qquad \qquad \qquad \qquad \qquad \qquad \qquad \qquad $$\\
$$+ \ [n]^{-2m+1} \ \sum_{j=2}^{m-1}\ \sum_{1\leq i_{1}\neq\ldots\neq i_{2m-j}\leq n}h(X_{i_{1}},\dots,X_{i_{j}},X_{i_{j+1}},\ldots,X_{i_{m}})\qquad \qquad \qquad  \qquad $$
$$\qquad \qquad \qquad \qquad \qquad \qquad \qquad \times\  h(X_{i_{1}},\dots,X_{i_{j}},X_{i_{m+1}},\ldots,X_{i_{2m-j}})$$\\
$$+\ [n]^{-2m+1}\ \sum_{1\leq i_{1}\neq\ldots\neq i_{2m-1}\leq n} h(X_{i_{1}},X_{i_{2}},\dots,X_{i_{m}}) \ h(X_{i_{1}},X_{i_{m+1}},\dots,X_{i_{2m-1}}). \qquad \qquad $$\\
\par
The first term and the second one, which obviously does not appear when $m=2$, in the latter equality will be seen to be negligible in probability (cf. Propositions 3 and 4), thus the third term becomes the main term that will play the main role in establishing Proposition 2.
\par
To complete the proof of Proposition 2 we shall state and prove the next three results, namely Propositions 3, 4 and Theorem 2.
\\ \\
\textbf{Proposition 3}. \emph{If} $\mathbb{E}|h(X_{1},\ldots,X_{m})|^{\frac{5}{3}}<\infty$, \emph{then, as} $n\rightarrow\infty$, \\ \\
$$[n]^{-2m+1} \sum_{1\leq i_{1}\neq\ldots\neq i_{m}\leq n}h^{2}(X_{i_{1}},\ldots,X_{i_{m}})\rightarrow 0\ \ \ a.s.$$ \\\\
\textit{\textbf{Proof of Proposition 3}}\\
From the fact that for $m\geq 2$, $\displaystyle{\frac{2m}{2m-1}}<\displaystyle{\frac{5}{3}}$, it follows that $$\mathbb{E}\left|h^{2}(X_{1},\ldots,X_{m})\right|^{\frac{m}{2m-1}}= \mathbb{E}\left|h(X_{1},\ldots,X_{m})\right|^{\frac{2m}{2m-1}}<\infty.$$
By this the proof of Proposition 3  follows from Theorem 1 of \cite{gin1}.
\\ \\
\textbf{Proposition 4}. For $m\geq 3$, \emph{If} $\mathbb{E}|h(X_{1},\ldots,X_{m})|^{\frac{5}{3}}<\infty$, \emph{then, as}  $n\rightarrow\infty$,\\
$$[n]^{-2m+1} \ \sum_{j=2}^{m-1}\ \sum_{1\leq i_{1}\neq\ldots\neq i_{2m-j}\leq n}h(X_{i_{1}},\dots,X_{i_{j}},X_{i_{j+1}},\ldots,X_{i_{m}})\qquad \qquad \qquad  \qquad $$
$$\qquad \qquad \qquad \qquad \qquad \qquad \qquad \times\  h(X_{i_{1}},\dots,X_{i_{j}},X_{i_{m+1}},\ldots,X_{i_{2m-j}})=o_{P}(1).$$\\
\textit{\textbf{Proof of Proposition 4}}\\
In order to prove Proposition 4 it suffices to show that as $n\rightarrow \infty$, for $j=2,\ldots,m-1,$ we have\\
\\ \\
$$[n]^{-2m+1} \ \sum_{1\leq i_{1}\neq\ldots\neq i_{2m-j}\leq n}h(X_{i_{1}},\dots,X_{i_{j}},X_{i_{j+1}},\ldots,X_{i_{m}})\qquad \qquad \qquad \qquad \qquad \qquad \qquad $$
$$\qquad \qquad \qquad \qquad \times\  h(X_{i_{1}},\dots,X_{i_{j}},X_{i_{m+1}},\ldots,X_{i_{2m-j}})=o_{P}(1).$$
Since the proof of the latter relation can be done by modifying, mutatis mutandis (cf. Appendix), that of the next theorem, i.e., Theorem 2, hence the detailed proof is given in Appendix.\\\\
\textbf{Theorem 2}. \emph{If} $\mathbb{E}|h(X_{1},\ldots,X_{m})|^{\frac{5}{3}}<\infty$ and $\tilde{h}_{1}(X_{1})\in DAN$ \emph{then, as} $n\rightarrow\infty$,\\
$$|\ [n]^{-2m+1}\ \sum_{1\leq i_{1}\neq\ldots\neq i_{2m-1}\leq n} h(X_{i_{1}},X_{i_{2}},\dots,X_{i_{m}})  h(X_{i_{1}},X_{i_{m+1}},\dots,X_{i_{2m-1}})$$
$$\qquad \qquad \qquad \qquad \qquad \qquad -\ \frac{1}{n}\sum_{i=1}^{n} \tilde{h}_{1}^{2}(X_{i})\ |=o_{P}(1). $$\\
\textit{\textbf{Proof of Theorem 2}}
\\
Before stating the proof of Theorem 2 we need the following definition and lemma which will play a crucial role in our proofs.
\\
\\
\textbf{Definition}. The Borel- measurable function $L(x_{1},\ldots,x_{m}):\mathbb{R}^{m}\rightarrow\mathbb{R}$, $m\geq 2$, with mean
$\mu=\mathbb{E}L(X_{1},\ldots,X_{m})$, is said to be \emph{degenerate}  if for every proper subset $\{\alpha_{1},\ldots,\alpha_{j}\}$ of $\{1,\ldots,m\}$, $j=1,\ldots,m-1$,  we have $$\mathbb{E}\textbf{(}L(X_{1},\ldots,X_{m})-\mu|X_{\alpha_{1}},\ldots,X_{\alpha_{j}}\textbf{)}=0 \ \ a.s.$$
We note in passing  that if $L$ were symmetric in its arguments, then the associated $U$-statistic with such a kernel would be a complete degenerate one. Hence our terminology for $L$ in this definition.
\\ \\
\textbf{Lemma 1}. \emph{If} $L:\mathbb{R}^{m}\rightarrow\mathbb{R}$, $m\geq 2$,  \emph{is degenerate  with mean} $\mu=\mathbb{E}L(X_{1},\ldots,X_{m})$

\emph{and} $\mathbb{E} L^{2}(X_{1},\ldots,X_{m})<\infty$, \emph{then,}
\\
$$\mathbb{E}\textbf{(}\ [n]^{-m} \sum_{1\leq i_{1}\neq \ldots \neq i_{m}\leq n}(L(X_{i_{1}},\ldots,X_{i_{m}})-\mu)\ \textbf{)}^{2}\leq \ [n]^{-m} \ \mathbb{E}\ \textbf{(}L(X_{1},\ldots,X_{m})-\mu \textbf{)}^{2}.$$
\\
\textit{\textbf{Proof of Lemma 1}}
\\ \\
Let $\hat{L}_{1\ldots m}:=\displaystyle{\frac{1}{m!}} \sum_{C_{m}}\ L_{\sigma_{1}\ldots \sigma_{m}}$, where $L_{\sigma_{1} \ldots \sigma_{m}}:= L(X_{\sigma_{1}},\ldots,X_{\sigma_{m}})$ and $C_{m}$ denotes the set of all permutations $\sigma_{1},\ldots,\sigma_{m}$ of $1,\ldots,m.$ It is clear that $$\sum_{1\leq i_{1}\neq \ldots \neq i_{m}\leq n}(\hat{L}_{i_{1},\ldots,i_{m}}-\mu)=\ \sum_{1\leq i_{1}\neq \ldots \neq i_{m}\leq n}(L_{i_{1},\ldots,i_{m}}-\mu).$$
Now write
\\ \\
$\mathbb{E}\ \textbf{(} \ [n]^{-m} \sum_{1\leq i_{1}\neq \ldots \neq i_{m}\leq n}(\hat{L}_{i_{1},\ldots,i_{m}}-\mu)\ \textbf{)}^{2}$
\\ \\ \\
$\textbf{=} \ ([n]^{-m})^{2} \sum_{1\leq i_{1}\neq \ldots \neq i_{m}\leq n} \mathbb{E}\textbf{(}\hat{L}_{i_{1},\ldots, i_{m}}-\mu\textbf{)}^{2}$
\\ \\ \\
$+\ ([n]^{-m})^{2} \sum_{j=1}^{m-1} \sum_{1\leq i_{1}\neq \ldots \neq i_{2m-j}\leq n} \mathbb{E}\textbf{[} \ \textbf{(}\hat{L}_{i_{1},\dots,i_{j},i_{j+1},\ldots,i_{m}}-\mu\textbf{)}$
$$\qquad \qquad \qquad \qquad \qquad \qquad \qquad \qquad \times \ \textbf{(}\hat{L}_{i_{1},\ldots,i_{j},i_{m+1},\ldots,i_{2m-j}}-\mu\textbf{)} \textbf{]}$$\\
$+\ ([n]^{-m})^{2}  \sum_{1\leq i_{1}\neq \ldots \neq i_{2m}\leq n} \mathbb{E}\ \textbf{(} \ (\hat{L}_{i_{1},\ldots,i_{m}}-\mu)\  (\hat{L}_{i_{m+1},\ldots,i_{2m}}-\mu)\ \textbf{)}$\\ \\
$= \ [n]^{-m} \ \mathbb{E}\ ( \hat{L}_{1,\ldots,m}-\mu )^{2}$\\
\\
$+ \ ([n]^{-m})^{2} \sum_{j=1}^{m-1} \sum_{1\leq i_{1}\neq \ldots \neq i_{2m-j}\leq n} \mathbb{E} \{ \mathbb{E}\textbf{[} \hat{L}_{i_{1},\dots,i_{j},i_{j+1},\ldots,i_{m}}-\mu\ |\ X_{i_{1}},\ldots,X_{i_{j}}\textbf{]}$
$$\qquad \qquad \qquad \qquad \qquad \qquad \qquad  \times \ \mathbb{E}\textbf{[}\hat{L}_{i_{1},\ldots,i_{j},i_{m+1},\ldots,i_{2m-j}}-\mu |\ X_{i_{1}},\ldots,X_{i_{j}} \textbf{]}\}$$\
$+\ ([n]^{-m})^{2}  \sum_{1\leq i_{1}\neq \ldots \neq i_{2m}\leq n} \mathbb{E}( \hat{L}_{i_{1},\ldots,i_{m}}-\mu)\  \mathbb{E}(\hat{L}_{i_{m+1},\ldots,i_{2m}}-\mu)$
\\ \\ \\
$=\ [n]^{-m} \ \mathbb{E}\ \textbf{(}\hat{L}_{1,\ldots,m}-\mu\textbf{)}^{2}$\\
\\
$\leq [n]^{-m}\  \mathbb{E}\ \textbf{(} L_{1,\ldots,m}-\mu\textbf{)}^{2}.$
\\ \\
The last inequality results from a well known inequality for sums of random variables followed by an application of Cauchy inequality provided that $\mathbb{E} L^{2}_{\sigma_{1},\ldots,\sigma_{m}}=\mathbb{E}L^{2}_{1\ldots m}$.
\par
It is easy to observe that when $L$ is symmetric in its arguments, the inequality in Lemma 1 becomes equality.
\\ \\
For further use in this proof, we consider the following setup:
\\ \\
$h_{1\ldots m}:=\ h(X_{1},\ldots,X_{m}),$
\\ \\
$h^{(m)}_{1\ldots m}:=\ h_{1\ldots m}\ \textbf{1}_{\displaystyle{(|h|\leq n^{\frac{3m}{5}})}},$
\\ \\
$h^{*}_{12\ldots 2m-1}:= \ h^{(m)}_{12\ldots m}\ h^{(m)}_{1 m+1 \ldots 2m-1},$
\\ \\
$\tilde{h}^{(m)}_{1}(x):=\ \mathbb{E}\ (\ h^{(m)}_{1\ldots m}|X_{1}=x \ ),$
\\ \\
$h^{(j)}_{1\ldots m}:=\ h^{(m)}_{1\ldots m } \ \textbf{1}_{\displaystyle{(|h^{(m)}| \leq n^{\frac{3j}{5}})}}$, $j=1,\ldots, m-1,$
\\ \\
$h^{(0)}_{1\ldots m}:=\ h^{(m)}_{1\ldots m}\ \textbf{1}_{\displaystyle{(|h^{(m)}| \leq \log(n))}},$
\\ \\
$h^{(\ell)}_{1\ldots m}:= \ h^{(m)}_{1\ldots m}\ \textbf{1}_{\displaystyle{(|\tilde{h}^{(m)}_{1}(x)| \leq n^{1/2}\ \ell(n))}},$
\\ \\
where $\textbf{1}_{A}$ denotes the indicator function of the set $A$ and $\ell(.)$ is a slowly varying function at infinity associated to $\tilde{h}_{1}(X_{1})$.
\par
In view of the above set up, observe that as $n\rightarrow \infty$
$$\mathbb{P}\ \textbf{(}\ \sum_{1\leq i_{1}\neq \ldots \neq i_{2m-1}\leq n}\ h_{i_{1} i_{2} \ldots i_{m} }\ h_{i_{1} i_{m+1}\ldots i_{2m-1}} \neq \sum_{1\leq i_{1}\neq \ldots \neq i_{2m-1}\leq n}\ h^{(m)}_{i_{1} i_{2} \ldots i_{m} }\ h^{(m)}_{i_{1} i_{m+1}\ldots i_{2m-1} } \textbf{)}$$
$$\leq \ n^{m}\ \mathbb{P}\ (\ |h_{1\ldots m}|>\ n^{\frac{3m}{5}}\ )\qquad \qquad \qquad \qquad \qquad \qquad \qquad \qquad \qquad \qquad \qquad $$
$$\leq \mathbb{E}\ [\ |h_{1 \ldots m}|^{\frac{5}{3}}\ 1_{(|h_{1 \ldots m}|>n^{\frac{3m}{5}})} \ ]\longrightarrow 0.\qquad \qquad \qquad \qquad \qquad \qquad \qquad \qquad \qquad \qquad \qquad \qquad $$
\\
Hence the asymptotic equivalency of the statistic of Theorem 2 and its truncated version i.e., $\sum_{1\leq i_{1}\neq \ldots \neq i_{2m-1}\leq n}\ h^{(m)}_{i_{1} i_{2} \ldots i_{m} }\ h^{(m)}_{i_{1} i_{m+1}\ldots i_{2m-1} }$  in probability.
\par
Having the asymptotic equivalency of the original statistic and its truncated version, to prove Theorem 2,  we shall proceed by working with the truncated version. Extending the idea of Hoffeding procedure to represent $U$-statistics in terms of complete \emph{degenerate} ones (cf. for example \cite{serf}), in our context in which due to lack of symmetry, our statistic of interest  i.e., $\sum_{1\leq i_{1}\neq \ldots \neq i_{2m-1}\leq n}\ h^{(m)}_{i_{1} i_{2} \ldots i_{m} }\ h^{(m)}_{i_{1} i_{m+1}\ldots i_{2m-1} }$ is not a $U$-statistic, by adding and subtracting required terms, we shall create a sequence of \emph{degenerate} statistics. Then by employing  proper new truncations and applying Lemma 1 we conclude the asymptotic negligibility of all of these degenerate statistics in probability (cf. Propositions 5, 6 and 7) except for the last group of them which are of the form of sums of i.i.d. random variables (cf. Remark 7). Among those the latter mentioned just one (cf. part (b) of Proposition 8) will  asymptotically in probability coincide $\displaystyle{\frac{1}{n}}\sum _{i=1}^{n}\tilde{h}^{2}_{1}(X_{i})$ and that will complete the proof of Theorem 2.
\par
Now by adding and subtracting required terms we  write
\\ \\
$$\sum_{1\leq i_{1}\neq \ldots \neq i_{2m-1}\leq n}\ h^{*}_{i_{1}\ldots i_{2m-1}}\qquad \qquad \qquad \qquad \qquad \qquad \qquad \qquad \qquad \qquad \qquad \qquad \qquad $$
$$\textbf{=}\ \sum_{1\leq i_{1}\neq \ldots \neq i_{2m-1}\leq n}\textbf{\{}\sum_{d=1}^{2m-1}(-1)^{2m-1-d} \sum_{1\leq j_{1}< \ldots < j_{d}\leq 2m-1} \mathbb{E}(h^{*}_{i_{1}\ldots i_{2m-1}}-\mathbb{E}(h^{*}_{i_{1}\ldots i_{2m-1}})|\ X_{i_{j_{1}}},\ldots X_{i_{j_{d}}})$$
$$\ \textbf{+}\sum_{c=1}^{2m-2}\sum_{1\leq k_{1}<\ldots<k_{c}\leq 2m-1}\sum_{d=1}^{c}(-1)^{c-d}\sum_{1\leq j_{1}< \ldots < j_{d}\leq c} \mathbb{E}(h^{*}_{i_{1}\ldots i_{2m-1}}-\mathbb{E}(h^{*}_{i_{1}\ldots i_{2m-1}})|\ X_{i_{k_{j_{1}}}},\ldots X_{i_{k_{j_{d}}}})$$
\\
$ $ $\ \  \textbf{+}\ \mathbb{E}\ (h^{*}_{i_{1} \ldots i_{2m-1} }) \textbf{\}}$
\\
$$\textbf{:=}\ \sum_{1\leq i_{1}\neq \ldots \neq i_{2m-1}\leq n} V(i_{1},\ldots,i_{2m-1})\ +\ \sum_{c=1}^{2m-2}\sum_{1\leq k_{1}<\ldots<k_{c}\leq 2m-1} \sum_{1\leq i_{1}\neq \ldots \neq i_{2m-1}\leq n} V(i_{k_{1}},\ldots,i_{k_{c}})$$
\\
$\ \ \ \textbf{+}\ \sum_{1\leq i_{1}\neq \ldots \neq i_{2m-1}\leq n} \mathbb{E}(h^{*}_{i_{1}\ldots i_{2m-1}}).$
\\ \\
\textbf{Proposition 5}. \emph{If} $\mathbb{E}\ |h_{1\ldots m}|^{\frac{5}{3}}<\infty$, \emph{then, as} $n\rightarrow \infty,$
$$[n]^{-2m+1}\sum_{1\leq i_{1}\neq \ldots \neq i_{2m-1}\leq n} V(i_{1},\ldots,i_{2m-1})\ =o_{P}(1).$$
\\
\textbf{\emph{Proof of Proposition 5}}
\\
For throughout use  $K$ will be a positive constant that may be different at each stage.
\par
 Since $V(i_{1},\ldots,i_{2m-1})$ posses the property of \emph{degeneracy} we can apply Lemma 1 for the associated statistics and write, for $\epsilon>0$,
\\ \\
$$\mathbb{P}\ (\ |\ [n]^{-2m+1}\sum_{1\leq i_{1}\neq \ldots \neq i_{2m-1}\leq n} V(i_{1},\ldots,i_{2m-1})|>\epsilon)\qquad $$
$$\leq \epsilon^{-2}\ \mathbb{E}\ \textbf{[}\ [n]^{-2m+1} \sum_{1\leq i_{1}\neq \ldots \neq i_{2m-1}\leq n} V(i_{1},\ldots,i_{2m-1})\ \textbf{]}^{2}  \qquad $$
$$\leq \epsilon^{-2}\ [n]^{-2m+1} \ \mathbb{E}\ [\ V(1,\ldots,2m-1) \ ]^{2}\qquad \qquad \qquad \qquad $$
$$\ \ \leq \ K \ \epsilon^{-2}\ [n]^{-2m+1}\ n^{2m-1} \ n^{-2m+1} \ \mathbb{E}\ [\ h^{(m)}_{12\ldots m}\ h^{(m)}_{1 m+1 \ldots 2m-1} \ ]^{2} $$
$$\leq \ K \ \epsilon^{-2}\ [n]^{-2m+1}\ n^{2m-1} \ n^{-2m+1} \ n^{\frac{7m}{5}} \ \mathbb{E}\ |\ h_{1 2\ldots m}\ |^{\frac{5}{3}}$$
$$\longrightarrow 0, \ \textrm{as} \ n\rightarrow \infty.\qquad \qquad \qquad \qquad \qquad \qquad \qquad $$
\\
\par
The estimation for $m\geq 3$ that occurs in our next proposition does not appear, and hence not needed, when $m=2$.
\\
\textbf{Proposition 6}. \emph{For} $m\geq 3$, \emph{if} $\mathbb{E}\ |\ h_{12\ldots m} \ |^{\frac{5}{3}}<\infty$, \emph{then as}, $n\rightarrow\infty$
\\
$$[n]^{-2m+1}\sum_{1\leq i_{1}\neq \ldots \neq i_{2m-1}\leq n} V(i_{k_{1}},\ldots,i_{k_{c}})\ =o_{P}(1),$$
\emph{where} $c=3,\ldots,2m-2$ \emph{and} $1\leq k_{1}<\ldots<k_{c}\leq 2m-1.$
\\
\\
\textbf{\emph{Proof of Proposition 6}}
\\
\\
Based on the way $i_{k_{1}},\ldots,i_{k_{c}}$ are distributed between
$h^{(m)}_{i_{1} i_{2} \ldots i_{m}}$ and $h^{(m)}_{i_{1} i_{m+1} \ldots i_{2m-1}}$ in two different cases when $k_{1}=1$ and $k_{1}\neq 1$, the proof is stated as follows.
\\ \\
\textbf{\emph{Case}} $k_{1}=1$
\\ \\
Let $s$ and $t$ be  respectively the number of elements of the sets $\{i_{k_{1}},\ldots, i_{k_{c}}\}\cap \{i_{1},i_{2},\ldots, i_{m}\} $ and $\{i_{k_{1}},\ldots, i_{k_{c}}\}\cap \{i_{1},i_{m+1},\ldots, i_{2m-1}\}$. It is clear that in this case, i.e., $k_{1}=1$, we have that $s,t\geq 1$ and $s+t=c+1$. Now define
\\
$$V^{T}(i_{k_{1}},\ldots,i_{k_{c}})= \sum_{d=1}^{c}(-1)^{c-d}\sum_{1\leq j_{1}<\ldots<j_{d}\leq c} \mathbb{E}(h^{*^{T}}_{i_{1} \ldots i_{2m-1}}-\mathbb{E}(h^{*^{T}}_{i_{1} \ldots i_{2m-1}})\ |\ x_{i_{k_{j_{1}}}},\ldots,x_{i_{k_{j_{d}}}}),\ (3)$$
\\
$$V^{T^{\prime}}(i_{k_{1}},\ldots,i_{k_{c}})= \sum_{d=1}^{c}(-1)^{c-d}\sum_{1\leq j_{1}<\ldots<j_{d}\leq c} \mathbb{E}(h^{*^{T^{\prime}}}_{i_{1} \ldots i_{2m-1}}-\mathbb{E}(h^{*^{T^{\prime}}}_{i_{1} \ldots i_{2m-1}})\ |\ x_{i_{k_{j_{1}}}},\ldots,x_{i_{k_{j_{d}}}}), \ (4)$$
\\
where $h^{*^{T}}_{i_{1} \ldots i_{2m-1}}=\ h^{(s)}_{i_{1} i_{2} \ldots i_{m}}\ h^{(m)}_{i_{1} i_{m+1} \ldots i_{2m-1}}$ and $h^{*^{T^{\prime}}}_{i_{1} \ldots i_{2m-1}}=\ h^{(s)}_{i_{1} i_{2} \ldots i_{m}}\ h^{(t)}_{i_{1} i_{m+1} \ldots i_{2m-1}}$.
\\
\par
Now observe that as $n\rightarrow\infty$
$$\mathbb{P} \ (\sum_{1\leq i_{1} \neq \ldots \neq i_{2m-1}\leq n}  V (i_{k_{1}},\ldots,i_{k_{c}})\ \neq  \sum_{1\leq i_{1} \neq \ldots \neq i_{2m-1}\leq n} V^{T^{\prime}}(i_{k_{1}},\ldots,i_{k_{c}})\ )$$
\\
$$\leq \mathbb{P} \ (\sum_{1\leq i_{1} \neq \ldots \neq i_{2m-1}\leq n}  V (i_{k_{1}},\ldots,i_{k_{c}})\ \neq  \sum_{1\leq i_{1} \neq \ldots \neq i_{2m-1}\leq n} V^{T}(i_{k_{1}},\ldots,i_{k_{c}})\ )$$
$$\qquad   +\ \mathbb{P} \ (\sum_{1\leq i_{1} \neq \ldots \neq i_{2m-1}\leq n}  V^{T} (i_{k_{1}},\ldots,i_{k_{c}})\ \neq  \sum_{1\leq i_{1} \neq \ldots \neq i_{2m-1}\leq n} V^{T^{\prime}}(i_{k_{1}},\ldots,i_{k_{c}})\ )$$
\\
$$\leq \ n^{s}\ \mathbb{P}\ (\ |h^{(m)}_{1 2 \ldots m}|>\ n^{\frac{3s}{5}} \ )+\ n^{t}\ \mathbb{P}\ (\ |h^{(m)}_{1 m+1 \ldots 2m-1}|>\ n^{\frac{3t}{5}} \ ) \qquad \qquad \ \ $$
\\
$$\leq \ \mathbb{E}\ [\ |h_{1 2\ldots m}|^{\frac{5}{3}}\ \textbf{1}_{(|h|>n^{\frac{3s}{5}})} \ ]+\ \mathbb{E}\ [\ |h_{1 m+1\ldots 2m-1}|^{\frac{5}{3}}\ \textbf{1}_{(|h|>n^{\frac{3t}{5}})} \ ]\longrightarrow 0.\qquad $$
\\
The latter relation suggests that $\sum_{1\leq i_{1} \neq \ldots \neq i_{2m-1}\leq n}  V (i_{k_{1}},\ldots,i_{k_{c}})$ and
\\
$\sum_{1\leq i_{1} \neq \ldots \neq i_{2m-1}\leq n}  V^{T^{\prime}} (i_{k_{1}},\ldots,i_{k_{c}})$ are asymptotically equivalent in probability.
\\
\par
Since $V^{T^{\prime}} (i_{k_{1}},\ldots,i_{k_{c}})$ is \emph{degenerate}, Markov inequality followed by  an application of Lemma 1 yields,
\\
$$\mathbb{P}\ (\ |\ [n]^{-2m+1}\sum_{1\leq i_{1}\neq \ldots \neq i_{2m-1}\leq n} V^{T^{\prime}}(i_{k_{1}},\ldots,i_{k_{c}})\ |\ >\epsilon \ ) \qquad \qquad $$
$$\leq \ \epsilon^{-2} \ \mathbb{E}\ [\ [n]^{-2m+1}\sum_{1\leq i_{1}\neq \ldots \neq i_{2m-1}\leq n} V^{T^{\prime}}(i_{k_{1}},\ldots,i_{k_{c}}) \ ]^{2}\qquad \qquad $$
$$\leq \ K \epsilon^{-2}\ [n-(2m-1-c)]^{-c}\ \mathbb{E}\ [\ h^{(s)}_{1 2 \ldots m}\ h^{(t)}_{1 m+1 \ldots 2m-1} \ ]^{2}\qquad \qquad $$
\\
$$\leq \ K \epsilon^{-2}\ [n-(2m-1-c)]^{-c}\ n^{c}\ n^{-c} \ n^{\frac{7(t+s)}{10}}\ \mathbb{E}\ |\ h_{1 2 \ldots m} \ |^{\frac{5}{3}}\qquad \qquad $$

$\ \ \ \ \longrightarrow 0,\ \textrm{as} \ n\rightarrow\infty.$
\\
The latter relation is true since when $c\geq 3$, we have  $-c+\displaystyle{\frac{7(t+s)}{10}}<0.$
\\
\\
\textbf{\emph{Case}} $k_{1}\neq 1$
\\
\\
Similarly to the previous case let $s$ and $t$ be  respectively the number of elements of the sets $\{i_{k_{1}},\ldots, i_{k_{c}}\}\cap \{i_{1},i_{2},\ldots, i_{m}\} $ and $\{i_{k_{1}},\ldots, i_{k_{c}}\}\cap \{i_{1},i_{m+1},\ldots, i_{2m-1}\}$. Clearly here  we have $s,t\geq 0$ and $s+t=c$. It is obvious that in this case $s,t$ can be zero but not simultaneously. More specifically, $(s=c,t=0)$ and $(s=0,t=c)$ can happen and due to their similarity we shall only treat $(s=c,t=0)$.
\\
\par
Let $V^{T}(i_{k_{1}},\ldots,i_{k_{c}})$ and $V^{T^{\prime}}(i_{k_{1}},\ldots,i_{k_{c}})$ be of the forms respectively (3) and (4), where $h^{*^{T}}_{i_{1} \ldots i_{2m-1}}=\ h^{(s)}_{i_{1} i_{2} \ldots i_{m}}\ h^{(m)}_{i_{1} i_{m+1} \ldots i_{2m-1}}$ and $h^{*^{T^{\prime}}}_{i_{1} \ldots i_{2m-1}}=\ h^{(s)}_{i_{1} i_{2} \ldots i_{m}}\ h^{(t)}_{i_{1} i_{m+1} \ldots i_{2m-1}}$.
\\
\par
Observe that as $n\rightarrow \infty$
\\
$$\mathbb{P}\ (\ \sum_{1\leq i_{1} \neq \ldots \neq i_{2m-1}\leq n}  V (i_{k_{1}},\ldots,i_{k_{c}})\ \neq \sum_{1\leq i_{1} \neq \ldots \neq i_{2m-1}\leq n}  V^{T^{\prime}} (i_{k_{1}},\ldots,i_{k_{c}}) \ )\qquad $$
\\
$$\leq \ \left\{
    \begin{array}{ll}
      n^{s}\ \mathbb{P}\ (\ |h^{(m)}_{1 2 \ldots m}|>\ n^{\frac{3s}{5}} \ )+\ n^{t}\ \mathbb{P}\ (\ |h^{(m)}_{1 m+1 \ldots 2m-1}|>\ n^{\frac{3t}{5}} \ ), & \hbox{s,t$>$ 0,\ s+t=c;} \\
      n^{c}\ \mathbb{P}(\ |h^{(m)}_{1 2 \ldots m}|>\ n^{\frac{3c}{5}} \ )+\ \mathbb{P}\ (\ |h^{(m)}_{1 m+1 \ldots 2m-1}|>\ \log(n) \ ), & \hbox{s=c,t=0} \\
          \end{array}
  \right.
$$
\\
$$\leq \left\{
         \begin{array}{ll}
           \mathbb{E}\ [\ |h_{1 2\ldots m}|^{\frac{5}{3}}\ \textbf{1}_{(|h|>n^{\frac{3s}{5}})} \ ]+\ \mathbb{E}\ [\ |h_{1 m+1\ldots 2m-1}|^{\frac{5}{3}}\ \textbf{1}_{(|h|>n^{\frac{3t}{5}})} \ ], & \hbox{s,t $>$ 0,\ s+t=c;} \\
          \mathbb{E}\ [\ |h_{1 2\ldots m}|^{\frac{5}{3}}\ \textbf{1}_{(|h|>n^{\frac{3c}{5}})} \ ]+\ \mathbb{P}\ (\ |h^{(m)}_{1 m+1 \ldots 2m-1}|>\ \log(n) \ ) , & \hbox{s=c,t=0}
         \end{array}
       \right.
$$
\\
$  \longrightarrow 0.$
\\ \\
Applying Markov inequality followed by an application Lemma 1 once again yields,
\\
$$\mathbb{P}\ (\ |\ [n]^{-2m+1} \sum_{1\leq i_{1}\neq \ldots \neq i_{2m-1}  \leq n}  \ V^{T^{\prime}}(i_{k_{1}},\ldots,i_{k_{c}}) |>\epsilon  )\qquad \qquad \qquad \qquad \qquad \qquad$$
\\
$$\leq K \epsilon^{-2}\ [n-(2m-1-c)]^{-c}\ n^{c}\ n^{-c}\ \mathbb{E}\ [\ h^{(s)}_{1 2 \ldots m}\ h^{(t)}_{1 m+1 \ldots 2m-1 }\ ]^{2}\qquad \qquad \qquad \ \ $$
\\
$$\leq \left\{
         \begin{array}{ll}
           K \epsilon^{-2}\ [n-(2m-1-c)]^{-c}\ n^{c}\ n^{-c} \ n^{\frac{7c}{10}}\ \mathbb{E}|h_{1 2 \ldots m }|^{\frac{5}{3}}, & \hbox{s,t$>$ 0, s+t=c;} \\
           K \epsilon^{-2}\ [n-(2m-1-c)]^{-c}\ n^{c}\ n^{-c}\ n^{\frac{7c}{10}}\ \log^{\frac{7}{6}}(n)\ \mathbb{E}|h_{1 2 \ldots m }|^{\frac{5}{3}} , & \hbox{s=c,t=0}
         \end{array}
       \right.
$$
$\longrightarrow 0,$ as $n\rightarrow\infty.$
\\
\\
This completes the proof of Proposition 6.
\\
\\
\textbf{Proposition 7}. \emph{If} $\mathbb{E}\ |h_{1 2 \ldots m}|^{\frac{5}{3}}<\infty$ \emph{and} $\tilde{h}_{1}(X_{1})\in DAN$, \emph{then, as} $n\rightarrow\infty$
\\
$$[n]^{-2m+1} \sum_{1\leq i_{1}\neq \ldots \neq i_{2m-1}  \leq n}  V(i_{k_{1}},i_{k_{2}})=o_{P}(1),$$
\\
\emph{where}, $1\leq k_{1}<k_{2}\leq 2m-1.$
\\\\
\textbf{\emph{Proof of Proposition 7}}
\\
As it was the case in the proof of the last proposition, we shall state the proof for two cases $k_{1}=1$ and $k_{1}\neq 1$ separately.
\\ \\
\textbf{\emph{Case $k_{1}=1$}
}\\
\\
Again let $s$ and $t$ be  respectively the number of elements of the sets $\{i_{k_{1}}, i_{k_{2}}\}\cap \{i_{1},i_{2},\ldots, i_{m}\} $ and $\{i_{k_{1}}, i_{k_{2}}\}\cap \{i_{1},i_{m+1},\ldots, i_{2m-1}\}$. It is clear that in this case we either have $(s=2,t=1)$ or $(s=1,t=2)$ which due to their similarity only $(s=2,t=1)$ will be treated as follows.
\par
 Define
$$V^{T}(i_{k_{1}},i_{k_{2}})= \sum_{d=1}^{2}(-1)^{2-d}\sum_{1\leq j_{1}<\ldots<j_{d}\leq 2} \mathbb{E}(h^{*^{T}}_{i_{1} \ldots i_{2m-1}}-\mathbb{E}(h^{*^{T}}_{i_{1} \ldots i_{2m-1}})\ |\ x_{i_{k_{j_{1}}}},\ldots,x_{i_{k_{j_{d}}}}),\ $$
\\
$$V^{T^{\prime}}(i_{k_{1}},i_{k_{2}})= \sum_{d=1}^{2}(-1)^{2-d}\sum_{1\leq j_{1}<\ldots<j_{d}\leq 2} \mathbb{E}(h^{*^{T^{\prime}}}_{i_{1} \ldots i_{2m-1}}-\mathbb{E}(h^{*^{T^{\prime}}}_{i_{1} \ldots i_{2m-1}})\ |\ x_{i_{k_{j_{1}}}},\ldots,x_{i_{k_{j_{d}}}}), \ $$
\\
where $h^{*^{T}}_{i_{1} \ldots i_{2m-1}}=\ h^{(2)}_{i_{1} i_{2} \ldots i_{m}}\ h^{(m)}_{i_{1} i_{m+1} \ldots i_{2m-1}}$ and $h^{*^{T^{\prime}}}_{i_{1} \ldots i_{2m-1}}=\ h^{(2)}_{i_{1} i_{2} \ldots i_{m}}\ h^{(\ell)}_{i_{1} i_{m+1} \ldots i_{2m-1}}$.
\\
\par
As $n\rightarrow\infty$, we have
\\
$$\mathbb{P}\ (\ \sum_{1\leq i_{1}\neq\ldots \neq i_{2m-1} \leq n} V(i_{k_{1}},i_{k_{2}})\neq \sum_{1\leq i_{1}\neq\ldots \neq i_{2m-1} \leq n} V^{T^{\prime}}(i_{k_{1}},i_{k_{2}}) \ )$$
$$\leq \ \mathbb{P}\ (\ \sum_{1\leq i_{1}\neq\ldots \neq i_{2m-1} \leq n} V(i_{k_{1}},i_{k_{2}})\neq \sum_{1\leq i_{1}\neq\ldots \neq i_{2m-1} \leq n} V^{T}(i_{k_{1}},i_{k_{2}}) \ ) $$
$$\qquad  +\ \mathbb{P}\ (\ \sum_{1\leq i_{1}\neq\ldots \neq i_{2m-1} \leq n} V^{T}(i_{k_{1}},i_{k_{2}})\neq \sum_{1\leq i_{1}\neq\ldots \neq i_{2m-1} \leq n} V^{T^{\prime}}(i_{k_{1}},i_{k_{2}}) \ )$$
\\
$$\leq \ n^{2}\ \mathbb{P}(\ |h^{(m)}_{12\ldots m}|>n^{6/5}  \ )+\ n\ \mathbb{P}(\ |\tilde{h}^{(m)}_{1}(X_{1})|> n^{1/2}\ \ell(n)  \ ) \ \ \ \ $$
\\
$$\leq \ \mathbb{E}\ [\ |h_{12\ldots m}|^{\frac{5}{3}} \ \textbf{1}_{(|h|>n^{6/5})}\ ] +\ n\ \mathbb{P}(\ |\tilde{h}^{(m)}_{1}(X_{1})|> n^{1/2}\ \ell(n)  \ ) \ \  $$
\\
$\ \ \ \ \ \ :=\ I_{1}(n)+\ I_{2}(n).$
\\ \\
It can be easily seen that as $n$ tends to infinity $I_{1}(n)\rightarrow 0$.\\

To deal with $I_{2}(n)$ we write
\\
$$n\ \mathbb{P}\ (\ |\tilde{h}^{(m)}_{1}(X_{1})|> n^{1/2}\ \ell(n)  \ )\qquad \qquad \qquad \qquad \qquad $$
$$\leq \ n\ \mathbb{P}\ (\ |\tilde{h}_{1}(X_{1})|> \frac{n^{1/2}\ \ell(n)}{2}  \ )\qquad \qquad \qquad \qquad \qquad $$
$$ \qquad \qquad +\ n\ \mathbb{P}\ ( \ |\ \mathbb{E}(h_{1 m+1 \ldots 2m-1}\ \textbf{1}_{(|h|>n^{\frac{3m}{5}})}\ |\ X_{1}) \ |> \frac{n^{1/2}\ \ell(n)}{2} )$$
$$\leq \ n\ \mathbb{P}\ (\ |\tilde{h}_{1}(X_{1})|> \frac{n^{1/2}\ \ell(n)}{2}  \ )\qquad \qquad \qquad \qquad \qquad $$
$$+\ 2\ n^{1/2}\ \ell^{-1}(n)\ \mathbb{E}\ [\ |\ h_{1 m+1 \ldots 2m-1} \ |\ \textbf{1}_{(|h|>n^{\frac{3m}{5}})}\   ]\ \ $$
$$\leq \ n\ \mathbb{P}\ (\ |\tilde{h}_{1}(X_{1})|> \frac{n^{1/2}\ \ell(n)}{2}  \ )\qquad \qquad \qquad \qquad \qquad $$
$$+\ 2\ n^{1/2}\ n^{-\frac{2m}{5}}\ \ell^{-1}(n) \ \mathbb{E} |\ h_{1 m+1 \ldots 2m-1} \ |^{\frac{5}{3}} \qquad \ \ \  $$
$$ \longrightarrow 0, \ \textrm{as}\ n\rightarrow\infty. \qquad \qquad \qquad \qquad \qquad \qquad \qquad \qquad
$$
The latter relation is true since $\tilde{h}_{1}(X_{1})\in DAN$ and $ m\geq 2$, and it means that $I_{2}(n)=o(1)$. Hence the asymptotic equivalency of $\sum_{1\leq i_{1}\neq\ldots \neq i_{2m-1} \leq n} V(i_{k_{1}},i_{k_{2}})$ and $ \sum_{1\leq i_{1} \neq \ldots \neq i_{2m-1} \leq n} V^{T^{\prime}}(i_{k_{1}},i_{k_{2}})$ in probability.
\\
\par
Before applying Lemma 1 for $[n]^{-2m+1}\ \sum_{1\leq i_{1} \neq \ldots \neq i_{2m-1} \leq n} V^{T^{\prime}}(i_{k_{1}},i_{k_{2}})$, since we know that $k_{1}=1$ and $s=2$, due to symmetry  of $h_{i_{1} i_{2} \ldots i_{m}}$, without loss of generality we assume that $k_{2}=2.$
\\
\par
Now for $\epsilon>0$, Markov inequality and Lemma 1 lead to
\\
$$\mathbb{P}(\ |\ [n]^{-2m+1}\ \sum_{1\leq i_{1}\neq \ldots \neq i_{2m-1} \leq n} V^{T^{\prime}}(i_{1},i_{2})\ |>\epsilon \ )\qquad \qquad \qquad \qquad \qquad \qquad \qquad $$
$$\leq\ K \ \epsilon^{-2}\ [n-(2m-3)]^{-2}\ n^{2}\ n^{-2}\ \mathbb{E}[\mathbb{E}(h^{(2)}_{1 2 \ldots m} h^{(\ell)}_{1 m+1 \ldots 2m-1}-\mathbb{E}(h^{(2)}_{1 2 \ldots m} h^{(\ell)}_{1 m+1 \ldots 2m-1}) | X_{1},X_{2})]^{2}$$
$$+ K \ \epsilon^{-2}\ [n-(2m-3)]^{-2}\ n^{2}\ n^{-2}\ \mathbb{E}[\mathbb{E}(h^{(2)}_{1 2 \ldots m} h^{(\ell)}_{1 m+1 \ldots 2m-1}-\mathbb{E}(h^{(2)}_{1 2 \ldots m} h^{(\ell)}_{1 m+1 \ldots 2m-1}) | X_{1})]^{2}$$
$$+ K \ \epsilon^{-2}\ [n-(2m-3)]^{-2}\ n^{2}\ n^{-2}\ \mathbb{E}[\mathbb{E}(h^{(2)}_{1 2 \ldots m} h^{(\ell)}_{1 m+1 \ldots 2m-1}-\mathbb{E}(h^{(2)}_{1 2 \ldots m} h^{(\ell)}_{1 m+1 \ldots 2m-1}) | X_{2})]^{2}$$
$$:=\ K \ \epsilon^{-2}\ [n-(2m-3)]^{-2}\ n^{2}\ J_{1}(n)\qquad \qquad \qquad \qquad \qquad \qquad \qquad \qquad \ \ $$
$$\ \ +\ K \ \epsilon^{-2}\ [n-(2m-3)]^{-2}\ n^{2}\ J_{2}(n)\qquad \qquad \qquad \qquad \qquad \qquad \qquad \qquad \ \ $$
$$\ \  +\ K \ \epsilon^{-2}\ [n-(2m-3)]^{-2}\ n^{2}\ J_{3}(n).\qquad \qquad \qquad \qquad \qquad \qquad \qquad \qquad \ \ $$
Considering that as $n\rightarrow\infty$, $[n-(2m-3)]^{-2}\ n^{2}\rightarrow 1$, we will show that $J_{1}(n),J_{2}(n),J_{3}(n)=o(1).$
\\
\par
To deal with $J_{1}(n)$ write
$$J_{1}(n)\leq \ n^{-2}\ \mathbb{E}[\ \mathbb{E}(\ h^{(2)}_{1 2 \ldots m} h^{(\ell)}_{1 m+1 \ldots 2m-1} \ |\ X_{1},X_{2} ) \ ]^{2} \qquad \qquad \qquad \qquad \qquad \qquad \qquad $$
$$=\ n^{-2}\ \mathbb{E}[\ \mathbb{E}^{2}(h^{(2)}_{12\ldots m}|\ X_{1},X_{2})\ \mathbb{E}^{2}(h^{(\ell)}_{1 m+1 \ldots 2m-1}|\ X_{1}) \ ]\qquad \qquad \qquad \ \ \ $$
$$\qquad  =\ n^{-2}\ \mathbb{E}[\ \mathbb{E}^{2}(h^{(2)}_{12\ldots m}|\ X_{1},X_{2})\ \mathbb{E}^{2}(h^{(m)}_{1 m+1 \ldots 2m-1}|\ X_{1})\ \textbf{1}_{(|\tilde{h}^{(m)}_{1}(X_{1})|\leq n^{1/2}\ \ell(n))} \ ]$$
$$\leq\ n^{-1}\ \ell^{2}(n)\ \mathbb{E}[\ h^{(2)}_{12\ldots m}  ]^{2} \qquad \qquad \qquad \qquad \qquad \qquad \qquad \qquad \qquad \ \  $$
$$\leq \ n^{-\frac{3}{5}}\ \ell^{2}(n)\ \mathbb{E}|\ h_{1 2 \ldots m} |^{\frac{5}{3}} \qquad \qquad \qquad \qquad \qquad \qquad \qquad \qquad \qquad \ $$

$\ \ \ \longrightarrow 0,\ \textrm{as} \ n\rightarrow\infty,$
\\
i.e., $J_{1}(n)=o(1).$ A similar argument yields, $J_{2}(n)=o(1)$, hence the details are omitted.
\par
As for $J_{3}(n)$ we write
\\
$$J_{3}(n)\leq\ n^{-2}\ \mathbb{E}[\ \mathbb{E}( h^{(2)}_{1 2 \ldots m} h^{(\ell)}_{1 m+1 \ldots 2m-1} \ |\ X_{2} ) \ ]^{2} \qquad \qquad \qquad \qquad \qquad \qquad \qquad \qquad $$
$$=\ n^{-2}\ \mathbb{E}\{\ \mathbb{E}[ \mathbb{E}(h^{(2)}_{1 2 \ldots m} h^{(\ell)}_{1 m+1 \ldots 2m-1}|X_{1},\ldots,X_{m} )\ |\ X_{2} \ ]  \}^{2} \qquad \ $$
$$=\ n^{-2}\ \mathbb{E}\{\ \mathbb{E}(h^{(2)}_{1 2 \ldots m}|\ X_{2})\ \mathbb{E}(h^{(\ell)}_{1 m+1 \ldots 2m-1}|\ X_{1})   \}^{2}\qquad  \qquad \qquad $$
$$\leq \ n^{-\frac{3}{5}}\ \ell^{2}(n)\ \mathbb{E}|\ h_{1 2 \ldots m} |^{\frac{5}{3}} \qquad \qquad \qquad \qquad \qquad \qquad \qquad \ \ \  \ $$

$\qquad \ \ \   \longrightarrow 0,\ \textrm{as} \ n\rightarrow\infty.$
\\
The latter relation means that $J_{3}(n)=o(1)$. By this the proof of Proposition 7 when $k_{1}=1$ is complete. \\

At this stage we give the proof of Proposition 7 when $k_{1}\neq 1.$
\\
\\
\textbf{\emph{Case $k_{1}\neq 1$}}
\\
\\
Once again let $s$ and $t$ be  respectively the number of elements of the sets $\{i_{k_{1}}, i_{k_{2}}\}\cap \{i_{1},i_{2},\ldots, i_{m}\} $ and $\{i_{k_{1}}, i_{k_{2}}\}\cap \{i_{1},i_{m+1},\ldots, i_{2m-1}\}$. It is obvious that in this case the possibilities are either $s=t=1$ or when $m\geq 3$, $(s=2,t=0)$ or $(s=0,t=2)$. We shall treat the cases $s=t=1$ and $(s=2,t=0)$ when $m\geq 3$, separately as follows.
\\
\\
\emph{Case $k_{1}\neq 1$: $s=t=1$  }
\\
\\
We note that here we have $k_{1}\in\{2,\ldots,m\}$ and $k_{2}\in\{m+1,\ldots,2m-1\}$.
\par
Now define
$$V^{T}(i_{k_{1}},i_{k_{2}})= \sum_{d=1}^{2}(-1)^{2-d}\sum_{1\leq j_{1}<\ldots<j_{d}\leq 2} \mathbb{E}(h^{*^{T}}_{i_{1} \ldots i_{2m-1}}-\mathbb{E}(h^{*^{T}}_{i_{1} \ldots i_{2m-1}})\ |\ x_{i_{k_{j_{1}}}},\ldots,x_{i_{k_{j_{d}}}}),\ $$
\\
$$V^{T^{\prime}}(i_{k_{1}},i_{k_{2}})= \sum_{d=1}^{2}(-1)^{2-d}\sum_{1\leq j_{1}<\ldots<j_{d}\leq 2} \mathbb{E}(h^{*^{T^{\prime}}}_{i_{1} \ldots i_{2m-1}}-\mathbb{E}(h^{*^{T^{\prime}}}_{i_{1} \ldots i_{2m-1}})\ |\ x_{i_{k_{j_{1}}}},\ldots,x_{i_{k_{j_{d}}}}), \ $$
\\
where $h^{*^{T}}_{i_{1} \ldots i_{2m-1}}=\ h^{(1)}_{i_{1} i_{2} \ldots i_{m}}\ h^{(m)}_{i_{1} i_{m+1} \ldots i_{2m-1}}$ and $h^{*^{T^{\prime}}}_{i_{1} \ldots i_{2m-1}}=\ h^{(1)}_{i_{1} i_{2} \ldots i_{m}}\ h^{(1)}_{i_{1} i_{m+1} \ldots i_{2m-1}}$. Now observe that as $n\rightarrow\infty$ we have
\\
$$\mathbb{P} \ (\sum_{1\leq i_{1} \neq \ldots \neq i_{2m-1}\leq n}  V (i_{k_{1}},i_{k_{2}})\ \neq  \sum_{1\leq i_{1} \neq \ldots \neq i_{2m-1}\leq n} V^{T^{\prime}}(i_{k_{1}},i_{k_{2}})\ )$$
\\
$$\leq \mathbb{P} \ (\sum_{1\leq i_{1} \neq \ldots \neq i_{2m-1}\leq n}  V (i_{k_{1}},i_{k_{2}})\ \neq  \sum_{1\leq i_{1} \neq \ldots \neq i_{2m-1}\leq n} V^{T}(i_{k_{1}},i_{k_{2}})\ )$$
$$\qquad   +\ \mathbb{P} \ (\sum_{1\leq i_{1} \neq \ldots \neq i_{2m-1}\leq n}  V^{T} (i_{k_{1}},i_{k_{2}})\ \neq  \sum_{1\leq i_{1} \neq \ldots \neq i_{2m-1}\leq n} V^{T^{\prime}}(i_{k_{1}},i_{k_{2}})\ )$$
\\
$$\leq \ 2 \ n\ \mathbb{P}\ (\ |h^{(m)}_{1 2 \ldots m}|>\ n^{3/5} \ )\qquad \qquad \qquad \qquad \qquad \qquad \qquad \qquad \qquad$$
$$\leq \ 2\ \mathbb{E}[\ |h_{12\ldots m}|^{\frac{5}{3}}\ \textbf{1}_{(|h|>n^{3/5})} \ ] \qquad \qquad \qquad \qquad \qquad \qquad \qquad \qquad \ \ \ \ \ $$

$\qquad \longrightarrow 0.$
\\
In view of the latter relation we apply Lemma 1 to $[n]^{-2m+1}\sum_{1\leq i_{1} \neq \ldots \neq i_{2m-1}\leq n} V^{T^{\prime}}(i_{k_{1}},i_{k_{2}})$ and we get
\\
$$\mathbb{P}(\ |[n]^{-2m+1}\sum_{1\leq i_{1} \neq \ldots \neq i_{2m-1}\leq n} V^{T^{\prime}}(i_{k_{1}},i_{k_{2}})|>\epsilon)\qquad \qquad \qquad \qquad \qquad \qquad \qquad $$
$$\leq \ K \ \epsilon^{-2}\ [n-(2m-3)]^{-2}\ n^{2}\ n^{-2}\ \mathbb{E}(h^{(1)}_{12\ldots m} h^{(1)}_{1 m+1 \ldots 2m-1})^{2}\qquad \qquad \qquad$$
$$\leq \ K \ \epsilon^{-2}\ [n-(2m-3)]^{-2}\ n^{2}\ n^{-2}\ n^{7/5} \mathbb{E}|h_{12\ldots m}|^{\frac{5}{3}}\qquad \qquad \qquad \qquad\ \qquad   $$

$\qquad \longrightarrow 0, \ \textrm{as} \ n\rightarrow\infty.$
\\
This completes the proof of Proposition 7 for the Case $k_{1}\neq 1$ when $s=t=1.$
\\
\\
\emph{Case $k_{1}\neq 1$:  ( $m\geq3$ ) $s=2,t=0$  \ }
\\
\\
In this case we first note that $k_{1},k_{2}\in\{2,\ldots,m\}$. Now define
$$V^{T}(i_{k_{1}},i_{k_{2}})= \sum_{d=1}^{2}(-1)^{2-d}\sum_{1\leq j_{1}<\ldots<j_{d}\leq 2} \mathbb{E}(h^{*^{T}}_{i_{1} \ldots i_{2m-1}}-\mathbb{E}(h^{*^{T}}_{i_{1} \ldots i_{2m-1}})\ |\ x_{i_{k_{j_{1}}}},\ldots,x_{i_{k_{j_{d}}}}),\ $$
\\
$$V^{T^{\prime}}(i_{k_{1}},i_{k_{2}})= \sum_{d=1}^{2}(-1)^{2-d}\sum_{1\leq j_{1}<\ldots<j_{d}\leq 2} \mathbb{E}(h^{*^{T^{\prime}}}_{i_{1} \ldots i_{2m-1}}-\mathbb{E}(h^{*^{T^{\prime}}}_{i_{1} \ldots i_{2m-1}})\ |\ x_{i_{k_{j_{1}}}},\ldots,x_{i_{k_{j_{d}}}}), \ $$
\\
where $h^{*^{T}}_{i_{1} \ldots i_{2m-1}}=\ h^{(2)}_{i_{1} i_{2} \ldots i_{m}}\ h^{(m)}_{i_{1} i_{m+1} \ldots i_{2m-1}}$ and $h^{*^{T^{\prime}}}_{i_{1} \ldots i_{2m-1}}=\ h^{(2)}_{i_{1} i_{2} \ldots i_{m}}\ h^{(0)}_{i_{1} i_{m+1} \ldots i_{2m-1}}$. Now observe that as $n\rightarrow\infty$
\\
$$\mathbb{P} \ (\sum_{1\leq i_{1} \neq \ldots \neq i_{2m-1}\leq n}  V(i_{k_{1}},i_{k_{2}})\ \neq  \sum_{1\leq i_{1} \neq \ldots \neq i_{2m-1}\leq n} V^{T^{\prime}}(i_{k_{1}},i_{k_{2}})\ )$$
\\
$$\leq \mathbb{P} \ (\sum_{1\leq i_{1} \neq \ldots \neq i_{2m-1}\leq n}  V (i_{k_{1}},i_{k_{2}})\ \neq  \sum_{1\leq i_{1} \neq \ldots \neq i_{2m-1}\leq n} V^{T}(i_{k_{1}},i_{k_{2}}\ )\ )$$
$$\qquad   +\ \mathbb{P} \ (\sum_{1\leq i_{1} \neq \ldots \neq i_{2m-1}\leq n}  V^{T} (i_{k_{1}},i_{k_{2}})\ \neq  \sum_{1\leq i_{1} \neq \ldots \neq i_{2m-1}\leq n} V^{T^{\prime}}(i_{k_{1}},i_{k_{2}})\ )$$
$$\leq \ n^{2}\ \mathbb{P}(\ |h^{(m)}_{12\ldots m}|>n^{6/5} \ )+\ \mathbb{P}(\ |h^{(m)}_{1 m+1 \ldots 2m-1}|>\log(n) \ )\qquad \ \ \ \ $$
$$\leq \ \mathbb{E}[\ |\ h_{12\ldots m} \ |^{\frac{5}{3}}\ \textbf{1}_{(|h|>n^{6/5})} \ ]+\ \mathbb{P}(\ |h_{1 m+1 \ldots 2m-1}|>\log(n) \ )\qquad $$

$\qquad \longrightarrow 0.$
\\
The latter relation together with \emph{degeneracy} of $V^{T^{\prime}}(i_{k_{1}},i_{k_{2}})$ enable us to use Lemma 1 once again  and arrive at
$$\mathbb{P}(|\ [n]^{-2m+1}\sum_{1\leq i_{1} \neq \ldots \neq i_{2m-1}\leq n} V^{T^{\prime}}(i_{k_{1}},i_{k_{2}}) \ |>\epsilon )\qquad \qquad \qquad $$
$$\leq \ K \ \epsilon^{-2} \ [n-(2m-3)]^{-2} \ n^{2}\ n^{-2}\ \mathbb{E}(h^{(2)}_{12\ldots m}\ h^{(0)}_{1 m+1 \ldots 2m-1})^{2}\ $$
$$\leq \ K \ \epsilon^{-2} \ [n-(2m-3)]^{-2} \ n^{2}\ n^{-\frac{3}{5}}\ \log^{7/6}(n)\ \mathbb{E}|h_{12\ldots m}|^{\frac{5}{3}} \qquad $$

$\ \qquad \longrightarrow 0, \ \textrm{as} \ n\rightarrow\infty.$
\\
Now the proof of Proposition 7 is complete.\\
\\
\textbf{Remark 7}. Before stating our next result we note in passing that when $k_{1}=1$ then
$[n]^{-2m+1}\sum_{1\leq i_{1}\neq \ldots \neq i_{2m-1}  \leq n}V(i_{k_{1}})$ is of the form
$$[n-(2m-2)]^{-1}\ \sum_{i_{1}\in\{1,\ldots,n\}/\{2,\ldots,2m-1\}}^{n} \mathbb{E}(h^{*}_{i_{1} 2 \ldots 2m-1}-\mathbb{E}(h^{*}_{i_{1} 2 \ldots 2m-1})\ |X_{i_{1}}), $$
otherwise, i.e., when for example $k_{1}=2 $ it has the following form
$$[n-(2m-2)]^{-1}\ \sum_{i_{2}\in\{1,\ldots,n\}/\{1,3,\ldots,2m-1\}}^{n} \mathbb{E}(h^{*}_{1 i_{2} 3 \ldots 2m-1}-\mathbb{E}(h^{*}_{1 i_{2} 3 \ldots 2m-1})\ |X_{i_{2}}),$$
and so on for $k_{1}\in\{2,\ldots,2m-1\}$.
\\
\\ \\ \\
\textbf{Proposition 8}. \emph{If} $\mathbb{E}|h_{1\ldots m}|^{\frac{5}{3}}<\infty$ \emph{and} $\tilde{h}_{1}(X_{1})\in DAN,$ \emph{then, as} $n\rightarrow \infty$
\\ \\
(a) $[n]^{-2m+1}\sum_{1\leq \neq i_{1}\neq \ldots \neq i_{2m-1} \leq n}V(i_{k_{1}})=o_{P}(1),\ for\ k_{1}\in\{2,\ldots,2m-1\},$
\\
$$(b)\ |\ [n-(2m-2)]^{-1}\ \sum_{i\in\{1,\ldots,n\}/\{2,\ldots,2m-1\}}^{n} \mathbb{E}(h^{*}_{i 2 \ldots 2m-1}-\mathbb{E}(h^{*}_{i 2 \ldots 2m-1})\ |X_{i}) \qquad $$
$$ \qquad  +\ \mathbb{E}(h^{*}_{12\ldots 2m-1})\ - \frac{1}{n}\ \sum_{i=1}^{n}\tilde{h}^{2}_{1}(X_{i})  \ |=o_{P}(1). \qquad \qquad \qquad \qquad \qquad \qquad  \qquad \qquad \qquad$$
\\
\textbf{\emph{Proof of Proposition 8}}
\\
\\
First we give the proof of part (a). Due to similarities, we state the proof only for $k_{1}=2$.
\par
Define
\\
$$V^{T}(i_{2})=\ \mathbb{E}(h^{*^{T}}_{i_{1} i_{2}\ldots i_{2m-1}}- \mathbb{E}(h^{*^{T}}_{i_{1} i_{2}\ldots i_{2m-1}})|\ X_{i_{2}} ),$$
$$V^{T^{\prime}}(i_{2})=\ \mathbb{E}(h^{*^{T^{\prime}}}_{i_{1} i_{2}\ldots i_{2m-1}}- \mathbb{E}(h^{*^{T^{\prime}}}_{i_{1} i_{2}\ldots i_{2m-1}})|\ X_{i_{2}} ),$$
where $h^{*^{T}}_{i_{1} i_{2}\ldots i_{2m-1}}=\ h^{(1)}_{i_{1} i_{2}\ldots i_{m}}\ h^{(m)}_{i_{1} i_{m+1}\ldots i_{2m-1}}$ and $h^{*^{T^{\prime}}}_{i_{1} i_{2}\ldots i_{2m-1}}=\ h^{(1)}_{i_{1} i_{2}\ldots i_{m}}\ h^{(0)}_{i_{1} i_{m+1}\ldots i_{2m-1}}.$
\\
Again observe that as $n\rightarrow \infty $
\\ \\
$$\mathbb{P}(\sum_{1\leq i_{1}\neq \ldots \neq i_{2m-1} \leq n}V(i_{2})\neq \sum_{1\leq i_{1}\neq \ldots \neq i_{2m-1} \leq n}V^{T^{\prime}}(i_{2})\ )\qquad \qquad \qquad \qquad $$
$$\leq\ \mathbb{P}(\sum_{1\leq i_{1}\neq \ldots \neq i_{2m-1} \leq n}V(i_{2})\neq \sum_{1\leq i_{1}\neq \ldots \neq i_{2m-1} \leq n}V^{T}(i_{2})\ )\qquad \qquad \qquad $$
$$+\ \mathbb{P}(\sum_{1\leq i_{1}\neq \ldots \neq i_{2m-1} \leq n} V^{T}(i_{2})\neq \sum_{1\leq i_{1}\neq \ldots \neq i_{2m-1} \leq n}V^{T^{\prime}}(i_{2})\ )  \qquad \qquad  $$
$$\leq\ n\ \mathbb{P}(|h^{(m)}_{12\ldots m}|>n^{3/5})+\ \mathbb{P}(|h^{(m)}_{1 m+1\ldots 2m-1}|>\log(n)) \qquad \qquad \qquad \ $$
$$\leq \ \mathbb{E}[\ |h_{1 2 \ldots m}|^{\frac{5}{3}} \ \textbf{1}_{(|h|>n^{3/5})} \ ]+\ \mathbb{P}(|h_{1 m+1\ldots 2m-1}|>\log(n)) \qquad \qquad \ \ $$

$\qquad \longrightarrow 0.$
\\
\\
An application of Markov inequality yields
$$\mathbb{P}\ (\ |\ [n]^{-2m+1}\sum_{1\leq i_{1}\neq \ldots \neq i_{2m-1} \leq n}V^{T^{\prime}}(i_{2})\ |>\epsilon )\qquad \qquad \qquad \qquad \qquad \qquad \qquad $$
$$\leq \ K \epsilon^{-2}\ [n-(2m-2)]^{-1}\ n \ n^{-1}\ \mathbb{E}(h^{(1)}_{12\ldots m}\ h^{(0)}_{1 m+1 \ldots 2m-1})^{2}\qquad \qquad \qquad \qquad $$
$$\leq \ K \epsilon^{-2}\ [n-(2m-2)]^{-1}\ n \ n^{-\frac{3}{10}}\ \log^{7/6}(n) \ \mathbb{E}|h_{1 2 \ldots m}|^{\frac{5}{3}}\qquad \qquad \qquad \qquad \ \ \  $$

$ \longrightarrow 0, \ \textrm{as} \ n\rightarrow \infty.$
\\ \\
This complete the proof of part (a).
\par
In the final stage of our proofs, to prove part (b) first define $\tilde{h^{*}}(x)=\ \mathbb{E}(h_{12\ldots m}\ \textbf{1}_{(|h|>n^{\frac{3m}{5}})}|X_{1}=x)$ and write
\\
$$|\frac{1}{n-2m+2} \ \sum_{i\in\{1,\ldots,n\}/\{2,\ldots,2m-1\}}^{n} \mathbb{E}(h^{*}_{i 2 \ldots 2m-1}-\mathbb{E}(h^{*}_{i 2 \ldots 2m-1})\ |X_{i})\qquad \qquad \qquad $$
$$ +\ \mathbb{E}(h^{*}_{12\ldots 2m-1})   - \ \frac{1}{n}\sum_{i=1}^{n}\tilde{h}^{2}_{1}(X_{i})\ |\qquad \qquad \qquad \qquad \qquad \qquad \qquad \qquad \qquad  $$
$$=|\frac{1}{n-2m+2} \ \sum_{i\in\{1,\ldots,n\}/\{2,\ldots,2m-1\}}^{n} \mathbb{E}(h^{*}_{i 2 \ldots 2m-1}\ |X_{i})\
 - \ \frac{1}{n}\sum_{i=1}^{n}\tilde{h}^{2}_{1}(X_{i})\ |\qquad \qquad \qquad $$
$$\leq |\frac{1}{n-2m+2} \ \sum_{i\in\{1,\ldots,n\}/\{2,\ldots,2m-1\}}^{n} \mathbb{E}(h^{*}_{i 2 \ldots 2m-1}\ |X_{i})\  \qquad \qquad \qquad \qquad \qquad \qquad $$
$$- \ \frac{1}{n-2m+2} \sum_{i=1}^{n}\tilde{h}^{2}_{1}(X_{i})\ |+\ \frac{2m-2}{n(n-2m+2)} \sum_{i=1}^{n}\tilde{h}^{2}_{1}(X_{i}) \qquad \qquad \qquad \qquad $$
$$\leq |\frac{1}{n-2m+2} \ \sum_{i\in\{1,\ldots,n\}/\{2,\ldots,2m-1\}}^{n} \mathbb{E}(h^{*}_{i 2 \ldots 2m-1}\ |X_{i}) \qquad \qquad \qquad \qquad \qquad \qquad $$
$$- \ \frac{1}{n-2m+2} \sum_{i\in\{1,\ldots,n\}/ \{2,\ldots,2m-1\}}\tilde{h}^{2}_{1}(X_{i})\ |+\frac{1}{n-2m+2} \sum_{i=2}^{2m-1}\tilde{h}^{2}_{1}(X_{i})\qquad $$
$$+\ \frac{2m-2}{n(n-2m+2)} \sum_{i=1}^{n}\tilde{h}^{2}_{1}(X_{i})\qquad \qquad \qquad \qquad \qquad \qquad \qquad \qquad \qquad \qquad  $$
$$=\ \frac{1}{n-2m+2}\ |\ \sum_{i\in\{1,\ldots,n\}/ \{2,\ldots,2m-1\}} [-\tilde{h^{*}}(X_{i})\ ]\ [2\tilde{h}^{(m)}_{1}(X_{i})+\ \tilde{h^{*}}(X_{i})] \  |\qquad \qquad \qquad $$
$$+\ \frac{1}{n-2m+2} \sum_{i=2}^{2m-1}\tilde{h}^{2}_{1}(X_{i})+\ \frac{2m-2}{n(n-2m+2)} \sum_{i=1}^{n}\tilde{h}^{2}_{1}(X_{i})\qquad \qquad \qquad \qquad  $$
$$\leq \ \frac{1}{n-2m+2}\ [\sum_{i\in\{1,\ldots,n\}/ \{2,\ldots,2m-1\}}\tilde{h}^{2}_{1}(X_{i})]^{1/2} [\sum_{i\in\{1,\ldots,n\}/ \{2,\ldots,2m-1\}}\tilde{h^{*}}^{2}(X_{i})]^{1/2}\qquad \qquad $$
$$+\ \frac{1}{n-2m+2}\sum_{i\in\{1,\ldots,n\}/ \{2,\ldots,2m-1\}}\tilde{h^{*}}^{2}(X_{i})\ +\ \frac{1}{n-2m+2} \sum_{i=2}^{2m-1}\tilde{h}^{2}_{1}(X_{i}) \qquad    $$
$$+\ \frac{2m-2}{n(n-2m+2)} \sum_{i=1}^{n}\tilde{h}^{2}_{1}(X_{i}).\qquad \qquad \qquad \qquad \qquad \qquad \qquad \qquad \qquad \ \ (5) $$
\\
It is easy to see that as $n\rightarrow\infty$, we have $\frac{1}{n-2m+2} \sum_{i=2}^{2m-1}\tilde{h}^{2}_{1}(X_{i})=o_{P}(1).$ Also in view of Corollary A, i.e., Raikov theorem, we have $\frac{2m-2}{n(n-2m+2)} \sum_{i=1}^{n}\tilde{h}^{2}_{1}(X_{i})=o_{P}(1)$, as $n\rightarrow\infty$. Hence, in view of (5), in order to complete the proof of part (b), it suffices to  show that as $n\rightarrow\infty$,
$$\frac{1}{n-2m+2}\sum_{i\in\{1,\ldots,n\}/ \{2,\ldots,2m-1\}}\tilde{h^{*}}^{2}(X_{i})=o_{P}(1).$$
To prove the latter relation we first use Markov inequality and conclude
\\
$$\mathbb{P}(\sum_{i\in\{1,\ldots,n\}/ \{2,\ldots,2m-1\}}\tilde{h^{*}}^{2}(X_{i})>\ \epsilon\ (n-2m+2) \ )\qquad \qquad \qquad \qquad $$
$$ \leq \ \epsilon^{-\frac{1}{2}}\ (n-2m+2)^{-\frac{1}{2}}\ \sum_{i\in\{1,\ldots,n\}/ \{2,\ldots,2m-1\}}\mathbb{E}\ |\ \tilde{h^{*}}^{2}(X_{i})\ |^{\frac{1}{2}}\qquad \qquad $$
$$\leq \ \epsilon^{-\frac{1}{2}}\ (n-2m+2)^{\frac{1}{2}} \ \mathbb{E}\ |\ \tilde{h^{*}}(X_{1})\ |\qquad \qquad \qquad \qquad \qquad \qquad \qquad $$
$$\leq \ \epsilon^{-\frac{1}{2}}\ (n-2m+2)^{\frac{1}{2}}\ n^{-\frac{1}{2}} \ n^{\frac{1}{2}}\  \mathbb{E}[\ |h_{12\ldots m}|\ \textbf{1}_{(|h|>n^{\frac{3m}{5}})} \ ]\qquad \qquad \ \ \ \ \  $$
$$\leq \ \epsilon^{-\frac{1}{2}}\ (n-2m+2)^{\frac{1}{2}}\ n^{-\frac{1}{2}}\ \mathbb{E}[\ |h_{12\ldots m}|^{\frac{5}{6m}+1}\ \textbf{1}_{(|h|>n^{\frac{3m}{5}})} \ ]\qquad \qquad \ \ \ \ \  $$

$\qquad \longrightarrow 0, \ \textrm{as} \ n\rightarrow\infty.$
\\
The latter relation is true since for $m\geq2$, we have that $\displaystyle{\frac{5}{6m}+1\leq \frac{5}{3}},$ and this completes the proof of part (b) and those of Proposition 8 and Theorem 2.
\\
\\
\textbf{Example}. Let $X_{1},X_{2}, \dots$, be a sequence of i.i.d. random variables with the density function
\\
$$ f(x)=\left\{
          \begin{array}{ll}
            |x-a|^{-3}, & \hbox{$|x-a|\geq1,\ a\neq 0,$} \\
            0\ \ \ \ \ \ \ \ \ \   , & \hbox{elsewhere.}
          \end{array}
        \right.$$
Consider the parameter $\theta=\mathbb{E}^{m}(X_{1})=a^{m}$, where $m\geq 1$ is a positive integer, and the kernel $h(X_{1},\ldots,X_{m})=\prod_{i=1}^{m} X_{i}$.
Then with $m, n$ satisfying $n\geq m $, the corresponding U-statistic is
$$U_{n}={n \choose m}^{-1} \sum_{C(n,m)} \prod_{j=1}^{m} X_{i_{j}}.$$
Simple calculation shows that $\tilde{h}_{1}(X_{1})=X_{1}\ a^{m-1}\ -\ a^{m}$.\\
\\
It is easy to check that $\mathbb{E}|h(X_{1},\dots,X_{m})|^{\frac{5}{3}}<\infty$ and that $\tilde{h}_{1}(X_{1})\in DAN$ (cf. Gut, \cite{gut}).
\par
For the pseudo-selfnormalized process
$$U_{[nt]}^{*}=\left\{
                    \begin{array}{ll}
                      \ \  \ 0,& \hbox{$0\leq t <\displaystyle{\frac{m}{n}}, $}\\ \\
                     {\frac {{ {[nt] \choose m}^{-1}\sum_{C([nt],m)} \prod_{j=1}^{m} X_{i_{j}}\ -\ a^{m}}}{  \textbf{(}\sum_{i=1}^{n}(X_{i} a^{m-1}\ -\ a^{m})^{2}\textbf{)}^{\frac{1}{2}}}}, & \hbox{$\displaystyle{\frac{m}{n}}\leq t\leq 1.$}\\
                    \end{array}
                  \right.$$
Nasari in \cite{nasa} concludes that $\frac{[nt]}{m} \ U_{[nt]}^{*} \longrightarrow_{d} W(t) \ \ \textrm{on} \ (D[0,1],\rho),$ where $\rho$ is the sup-norm for functions in $D[0,1]$ and $\{W(t),0\leq t\leq 1\}$ is a standard Wiener process.
\par
The studentized U-process based on $U_{n}$ here is defined as follows.
$$\displaystyle{U}_{[nt]}^{\emph{stu}}=\left\{
\begin{array}{ll}
\ 0,& \hbox{$0\leq t <\displaystyle{\frac{m}{n}}, $}\\
                    \displaystyle{ \frac{{[nt] \choose m}^{-1} \sum_{C([nt],m)} \prod_{j=1}^{m} X_{i_{j}}-\theta} {\sqrt{(n-1)\sum_{i=1}^{n}(U^{i}_{n-1}-U_{n})^{2}}}}, & \hbox{$\displaystyle{\frac{m}{n}}\leq t\leq 1,$}
                     \end{array}
                     \right.$$\\
  where, by (2),
  $$(n-1)\ \sum_{i=1}^{n}(U^{i}_{n-1}-\ U_{n})^{2}
=\frac{m^{2}(n-1)}{(n-m)^{2}}\textbf{\{}\sum_{i=1}^{n}\  X^{2}_{i} [{n-1 \choose m-1}^{-1}\sum_{\substack{1\leq i_{2}<\ldots<i_{m}\leq n \\ i_{2},\ldots,i_{m}\neq i}} \prod_{j=2}^{m}\  X_{i_{j}}\ ]^{2}$$
$$ \qquad -\ n \ [{n \choose m}^{-1} \sum_{C(n,m)} \prod_{j=1}^{m} X_{i_{j}}]^{2}\textbf{\}}.$$
In view of $\displaystyle{U}_{[nt]}^{\emph{stu}}$ and $U_{[nt]}^{*}$, our Main Theorem is applicable for $\displaystyle{U}_{[nt]}^{\emph{stu}}$ provided Theorem 1 continues hold true in this case. Hence, part (c) of Main Theorem implies that
$\displaystyle{[nt]\ U^{\emph{stu}}_{[nt]}} \rightarrow_{d}\ W(t)$ on ($D$[0,1],$\rho$), where $\rho$ is the sup-norm for functions in $D[0,1]$ and $\{W(t),0\leq t\leq 1\}$ is a standard Wiener process.
\newpage
\appendix{\textbf{Appendix: Proof of Proposition 4}}
\\
\\
As it was mentioned before, the proof of this proposition can be done by modifying that of Theorem 2, except for that some of the steps are not required. This is due to the presence  of the extra term of $n$ with negative power i.e.,  $n^{-j+1}$  in this proposition,  where  $j=2,\ldots,
 m-1$, and $m\geq 3$. It is clear that among the statistics in proposition 4 the one associated to $j=2$ has the largest extra term of  $n^{-1}$. Hence, we shall only show that as $n\rightarrow \infty$,
$$[n]^{-2m+1} \ \sum_{1\leq i_{1}\neq\ldots\neq i_{2m-2}\leq n}h(X_{i_{1}},X_{i_{2}},X_{i_{3}},\ldots,X_{i_{m}})\qquad \qquad \qquad \qquad \qquad \qquad \qquad $$
$$\qquad \qquad \qquad \qquad \times\  h(X_{i_{1}},X_{i_{2}},X_{i_{m+1}},\ldots,X_{i_{2m-2}})=o_{P}(1).\qquad (\textrm{I}) $$
\par
To prove (I), consider the following setup:
\\ \\
$h_{1\ldots m}:=\ h(X_{1},\ldots,X_{m}),$
\\ \\
$h^{(m)}_{1\ldots m}:=\ h_{1\ldots m}\ \textbf{1}_{\displaystyle{(|h|\leq n^{\frac{3m}{5}})}},$
\\ \\
$h^{**}_{12\ldots 2m-2}:= \ h^{(m)}_{12 3\ldots m}\ h^{(m)}_{1  2 \ m+1 \ldots 2m-2},$
\\ \\
$h^{(j)}_{1\ldots m}:=\ h^{(m)}_{1\ldots m } \ \textbf{1}_{\displaystyle{(|h^{(m)}| \leq n^{\frac{3j}{5}})}}$, $j=1,\ldots, m-1,$
\\ \\
$h^{(0)}_{1\ldots m}:=\ h^{(m)}_{1\ldots m}\ \textbf{1}_{\displaystyle{(|h^{(m)}| \leq \log(n))}},$
\\ \\
where $\textbf{1}_{A}$ is the indicator function of the set $A$.
\par
Now observe that as $n\rightarrow\infty$
$$\mathbb{P}\ \textbf{(}\ \sum_{1\leq i_{1}\neq \ldots \neq i_{2m-2}\leq n}\ h_{i_{1} i_{2} i_{3} \ldots i_{m} }\ h_{i_{1} i_{2} i_{m+1}\ldots i_{2m-2}} \neq \sum_{1\leq i_{1}\neq \ldots \neq i_{2m-2}\leq n}\ h^{(m)}_{i_{1} i_{2} i_{3} \ldots i_{m} }\ h^{(m)}_{i_{1} i_{2} i_{m+1}\ldots i_{2m-2} } \textbf{)}$$
$$\leq \ n^{m}\ \mathbb{P}\ (\ |h_{1\ldots m}|>\ n^{\frac{3m}{5}}\ ) \qquad \qquad \qquad \qquad \qquad \qquad \qquad \qquad \qquad \qquad \qquad \qquad $$
$$\leq \mathbb{E}\ [\ |h_{1 \ldots m}|^{\frac{5}{3}}\ 1_{(|h_{1 \ldots m}|>n^{\frac{3m}{5}})} \ ]\longrightarrow 0. \qquad \qquad \qquad \qquad \qquad \qquad \qquad \qquad \qquad \qquad$$
\\
In view of the latter asymptotic equivalency and our setup, in order to prove (I), we need to show that as $n\rightarrow \infty$,
$$[n]^{-2m+1} \ \sum_{1\leq i_{1}\neq\ldots\neq i_{2m-2}\leq n}h^{**}_{i_{1},i_{2},\ldots,i_{2m-2}}=o_{P}(1).$$
\par
Similarly to what we had in the proof of Theorem 2 we write
$$\sum_{1\leq i_{1}\neq \ldots \neq i_{2m-2}\leq n}\ h^{**}_{i_{1}\ldots i_{2m-2}}\qquad \qquad \qquad \qquad \qquad \qquad \qquad \qquad \qquad \qquad \qquad \qquad \qquad $$
$$\textbf{=}\ \sum_{1\leq i_{1}\neq \ldots \neq i_{2m-2}\leq n}\textbf{\{}\sum_{d=1}^{2m-2}(-1)^{2m-2-d} \sum_{1\leq j_{1}< \ldots < j_{d}\leq 2m-2} \mathbb{E}(h^{**}_{i_{1}\ldots i_{2m-2}}-\mathbb{E}(h^{**}_{i_{1}\ldots i_{2m-2}})|\ X_{i_{j_{1}}},\ldots X_{i_{j_{d}}})$$
$$\ \textbf{+}\sum_{c=1}^{2m-3}\sum_{1\leq k_{1}<\ldots<k_{c}\leq 2m-2}\sum_{d=1}^{c}(-1)^{c-d}\sum_{1\leq j_{1}< \ldots < j_{d}\leq c} \mathbb{E}(h^{**}_{i_{1}\ldots i_{2m-1}}-\mathbb{E}(h^{**}_{i_{1}\ldots i_{2m-2}})|\ X_{i_{k_{j_{1}}}},\ldots X_{i_{k_{j_{d}}}})$$
\\
$ $ $\ \  \textbf{+}\ \mathbb{E}\ (h^{**}_{i_{1} \ldots i_{2m-2} }) \textbf{\}}$
\\
$$\textbf{:=}\ \sum_{1\leq i_{1}\neq \ldots \neq i_{2m-2}\leq n} V^{*}(i_{1},\ldots,i_{2m-2})\ +\ \sum_{c=1}^{2m-3}\sum_{1\leq k_{1}<\ldots<k_{c}\leq 2m-2} \sum_{1\leq i_{1}\neq \ldots \neq i_{2m-2}\leq n} V^{*}(i_{k_{1}},\ldots,i_{k_{c}})$$
\\
$\ \ \ \textbf{+}\ \sum_{1\leq i_{1}\neq \ldots \neq i_{2m-2}\leq n} \mathbb{E}(h^{**}_{i_{1}\ldots i_{2m-2}}).$
\\ \\
To prove (I), we shall show the asymptotic negligibility of all of the above terms in the next three propositions.
\\
\textbf{Proposition 4.1}: \emph{If} $\mathbb{E}\ |h_{1\ldots m}|^{\frac{5}{3}}<\infty$, \emph{then, as} $n\rightarrow \infty$
$$[n]^{-2m+1}\sum_{1\leq i_{1}\neq \ldots \neq i_{2m-2}\leq n} V^{*}(i_{1},\ldots,i_{2m-2})\ =o_{P}(1).$$
\\
\textbf{\emph{Proof of Proposition 4.1}}
\\
For throughout use  $K$ will be a positive constant that may be different at each stage.
\par
 Since $V^{*}(i_{1},\ldots,i_{2m-2})$ posses the property of \emph{degeneracy} we can apply Lemma 1 following a Markov inequality for the associated statistic and write, for $\epsilon>0$,
\\ \\
$$\mathbb{P}\ (\ |\ [n]^{-2m+2}\sum_{1\leq i_{1}\neq \ldots \neq i_{2m-2}\leq n} V^{*}(i_{1},\ldots,i_{2m-2})|>\epsilon\ (n-2m+2)\ )\qquad $$
$$\leq \epsilon^{-2}\ \ (n-2m+2)^{-2} \ \mathbb{E}\ \textbf{[}\ [n]^{-2m+2} \sum_{1\leq i_{1}\neq \ldots \neq i_{2m-1}\leq n} V^{*}(i_{1},\ldots,i_{2m-2})\ \textbf{]}^{2}  \qquad $$
$$\leq \epsilon^{-2}\  \ (n-2m+2)^{-2}\ [n]^{-2m+2} \ \mathbb{E}\ [\ V^{*}(1,\ldots,2m-2) \ ]^{2}\qquad \qquad \qquad \qquad $$
$$\leq \ K \ \epsilon^{-2}\ (n-2m+2)^{-2} \ [n]^{-2m+2}\ n^{2m} \ n^{-2m} \ \mathbb{E}\ [\ h^{(m)}_{12\ 3\ldots m}\ h^{(m)}_{1 2\ m+1 \ldots 2m-2} \ ]^{2} \qquad \qquad \qquad $$
$$\leq \ K \ \epsilon^{-2}\ (n-2m+2)^{-2} \ [n]^{-2m+2}\ n^{2m} \ n^{-2m} \ n^{\frac{7m}{5}} \ \mathbb{E}\ |\ h_{1 2\ldots m}\ |^{\frac{5}{3}} \qquad \qquad \qquad \qquad \qquad $$
$$\longrightarrow 0, \ \textrm{as} \ n\rightarrow \infty.\qquad \qquad \qquad \qquad \qquad \qquad \qquad \qquad \qquad \qquad \qquad \qquad$$
\\
\\
\textbf{Proposition 4.2}. \emph{If} $\mathbb{E}\ |h_{1\ldots m}|^{\frac{5}{3}}<\infty$, \emph{then, as} $n\rightarrow \infty,$
$$[n]^{-2m+1}\sum_{1\leq i_{1}\neq \ldots \neq i_{2m-2}\leq n} V(i_{k_{1}},\ldots,i_{k_{c}})\ =o_{P}(1),$$
\emph{where} $c=2,\ldots,2m-3$ \emph{and} $1\leq k_{1}<\ldots<k_{c}\leq 2m-2.$
\\
\\
\textbf{\emph{Proof of Proposition 4.2}}
\\
The proof will be stated in three cases according to the values of $k_{1}$ and $k_{2}$ as follows.
\\
\par
\emph{Case $k_{1}=1$ and $k_{2}=2$}
\\
Let $s$ and $t$ be  respectively the number of elements of the sets $\{i_{k_{1}},\ldots, i_{k_{c}}\}\cap \{i_{1},i_{2},i_{3},\ldots, i_{m}\} $ and $\{i_{k_{1}},\ldots, i_{k_{c}}\}\cap \{i_{1},i_{2},i_{m+1},\ldots, i_{2m-1}\}$. It is clear that in this case, i.e., $k_{1}=1$ and $k_{2}=2$, we have that $s,t\geq 2$ and $s+t=c+2$. Now define
\\
$$V^{*^{T}}(i_{k_{1}},\ldots,i_{k_{c}})= \sum_{d=1}^{c}(-1)^{c-d}\sum_{1\leq j_{1}<\ldots<j_{d}\leq c} \mathbb{E}(h^{**^{T}}_{i_{1} \ldots i_{2m-2}}-\mathbb{E}(h^{**^{T}}_{i_{1} \ldots i_{2m-2}})\ |\ x_{i_{k_{j_{1}}}},\ldots,x_{i_{k_{j_{d}}}}),\ $$
\\
$$V^{*T^{\prime}}(i_{k_{1}},\ldots,i_{k_{c}})= \sum_{d=1}^{c}(-1)^{c-d}\sum_{1\leq j_{1}<\ldots<j_{d}\leq c} \mathbb{E}(h^{**^{T^{\prime}}}_{i_{1} \ldots i_{2m-2}}-\mathbb{E}(h^{**^{T^{\prime}}}_{i_{1} \ldots i_{2m-2}})\ |\ x_{i_{k_{j_{1}}}},\ldots,x_{i_{k_{j_{d}}}}), \ $$
\\
where $h^{**^{T}}_{i_{1}  \ldots i_{2m-2}}=\ h^{(s)}_{i_{1} i_{2}\ i_{3} \ldots i_{m}}\ h^{(m)}_{i_{1}  i_{2}\  i_{m+1} \ldots i_{2m-2}}$ and $h^{**^{T^{\prime}}}_{i_{1} \ldots i_{2m-2}}=\ h^{(s)}_{i_{1} i_{2} \ i_{3} \ldots i_{m}}\ h^{(t)}_{i_{1} i_{2}\ i_{m+1} \ldots i_{2m-2}}$.
\\
\par
Now observe that as $n\rightarrow\infty$
$$\mathbb{P} \ (\sum_{1\leq i_{1} \neq \ldots \neq i_{2m-2}\leq n}  V^{*} (i_{k_{1}},\ldots,i_{k_{c}})\ \neq  \sum_{1\leq i_{1} \neq \ldots \neq i_{2m-2}\leq n} V^{*^{T^{\prime}}}(i_{k_{1}},\ldots,i_{k_{c}})\ )$$
\\
$$\leq \mathbb{P} \ (\sum_{1\leq i_{1} \neq \ldots \neq i_{2m-2}\leq n}  V^{*} (i_{k_{1}},\ldots,i_{k_{c}})\ \neq  \sum_{1\leq i_{1} \neq \ldots \neq i_{2m-2}\leq n} V^{*^{T}}(i_{k_{1}},\ldots,i_{k_{c}})\ )$$
$$\qquad   +\ \mathbb{P} \ (\sum_{1\leq i_{1} \neq \ldots \neq i_{2m-2}\leq n}  V^{*^{T}} (i_{k_{1}},\ldots,i_{k_{c}})\ \neq  \sum_{1\leq i_{1} \neq \ldots \neq i_{2m-2}\leq n} V^{*^{T^{\prime}}}(i_{k_{1}},\ldots,i_{k_{c}})\ )$$
\\
$$\leq \ n^{s}\ \mathbb{P}\ (\ |h^{(m)}_{1 2 \ 3 \ldots m}|>\ n^{\frac{3s}{5}} \ )+\ n^{t}\ \mathbb{P}\ (\ |h^{(m)}_{1 2\ m+1 \ldots 2m-1}|>\ n^{\frac{3t}{5}} \ ) \qquad \qquad \ \ $$
\\
$$\leq \ \mathbb{E}\ [\ |h_{1 2\ 3\ldots m}|^{\frac{5}{3}}\ \textbf{1}_{(|h|>n^{\frac{3s}{5}})} \ ]+\ \mathbb{E}\ [\ |h_{1 2\ m+1\ldots 2m-2}|^{\frac{5}{3}}\ \textbf{1}_{(|h|>n^{\frac{3t}{5}})} \ ]\longrightarrow 0.\qquad $$
\\
The latter relation suggests that $\sum_{1\leq i_{1} \neq \ldots \neq i_{2m-2}\leq n}  V^{*} (i_{k_{1}},\ldots,i_{k_{c}})$ and
\\
$\sum_{1\leq i_{1} \neq \ldots \neq i_{2m-2}\leq n}  V^{*^{T^{\prime}}} (i_{k_{1}},\ldots,i_{k_{c}})$ are asymptotically equivalent in probability.
\\
\par
Since $V^{*^{T^{\prime}}} (i_{k_{1}},\ldots,i_{k_{c}})$ is \emph{degenerate}, Markov inequality followed by  an application of Lemma 1 yields,
\\
$$\mathbb{P}\ (\ |\ [n]^{-2m+2}\sum_{1\leq i_{1}\neq \ldots \neq i_{2m-2}\leq n} V^{*^{T^{\prime}}}(i_{k_{1}},\ldots,i_{k_{c}})\ |\ >\epsilon\ (n-2m+2) \ ) \qquad \qquad $$
$$\leq \ \epsilon^{-2} \ (n-2m+2)^{-2}\ \mathbb{E}\ [\ [n]^{-2m+2}\sum_{1\leq i_{1}\neq \ldots \neq i_{2m-2}\leq n} V^{*^{T^{\prime}}}(i_{k_{1}},\ldots,i_{k_{c}}) \ ]^{2}\qquad \qquad $$
$$\leq \ K \epsilon^{-2}\ (n-2m+2)^{-2} \ [n-(2m-2-c)]^{-c}\ \mathbb{E}\ [\ h^{(s)}_{1 2 3\ldots m}\ h^{(t)}_{1 2\  m+1 \ldots 2m-1} \ ]^{2}\qquad \qquad $$
\\
$$\leq \ K \epsilon^{-2}\ (n-2m+2)^{-2}  \ [n-(2m-1-c)]^{-c}\ n^{c+2}\ n^{-c-2} \ n^{\frac{7(t+s)}{10}}\ \mathbb{E}\ |\ h_{1 \ldots m} \ |^{\frac{5}{3}}\qquad \qquad $$

$\ \ \ \ \longrightarrow 0,\ \textrm{as} \ n\rightarrow\infty.$
\\
This completes the proof of Proposition 4.2 when $k_{1}=1$ and $k_{2}=2$.
\\ \\

\emph{Case either $k_{1}\neq 1$ or $k_{2}\neq 2$}
\\
\\
Let $s,t$ be as were defined in the previous case and note that here we have that  $s,t\geq 1$ and $s+t=c+1$. The proof of Proposition 4.2 in this case results from a similar argument to what   was given for the previous case, hence the details are omitted.
\\ \\

\emph{Case $k_{1}\neq 1,k_{2}\neq 2$ }
\\
\\
Let  $s,t$ be as what were defined in the previous two cases and note that in this case we have $s,t\geq 0$ and $s+t=c$. Also let $V^{*^{T}}$ and $V^{*^{T^{\prime }}}$ as they were defined in the case $k_{1}=1,k_{2}=2$ and observe that as $n\rightarrow\infty$ we have
$$\mathbb{P} \ (\sum_{1\leq i_{1} \neq \ldots \neq i_{2m-2}\leq n}  V^{*} (i_{k_{1}},\ldots,i_{k_{c}})\ \neq  \sum_{1\leq i_{1} \neq \ldots \neq i_{2m-2}\leq n} V^{*^{T^{\prime}}}(i_{k_{1}},\ldots,i_{k_{c}})\ )$$
\\
$$\leq \mathbb{P} \ (\sum_{1\leq i_{1} \neq \ldots \neq i_{2m-2}\leq n}  V^{*} (i_{k_{1}},\ldots,i_{k_{c}})\ \neq  \sum_{1\leq i_{1} \neq \ldots \neq i_{2m-2}\leq n} V^{*^{T}}(i_{k_{1}},\ldots,i_{k_{c}})\ )$$
$$\qquad   +\ \mathbb{P} \ (\sum_{1\leq i_{1} \neq \ldots \neq i_{2m-2}\leq n}  V^{*^{T}} (i_{k_{1}},\ldots,i_{k_{c}})\ \neq  \sum_{1\leq i_{1} \neq \ldots \neq i_{2m-2}\leq n} V^{*^{T^{\prime}}}(i_{k_{1}},\ldots,i_{k_{c}})\ )$$
\\
$$\leq \ \left\{
    \begin{array}{ll}
      n^{s}\ \mathbb{P}\ (\ |h^{(m)}_{1 2 \ 3 \ldots m}|>\ n^{\frac{3s}{5}} \ )+\ n^{t}\ \mathbb{P}\ (\ |h^{(m)}_{1 2\ m+1 \ldots 2m-2}|>\ n^{\frac{3t}{5}} \ ), & \hbox{s,t $>$ 0,\ s+t=c;} \\
      n^{c}\ \mathbb{P}(\ |h^{(m)}_{1 2 \ 3 \ldots m}|>\ n^{\frac{3c}{5}} \ )+\ \mathbb{P}\ (\ |h^{(m)}_{1 2\ m+1 \ldots 2m-2}|>\ \log(n) \ ), & \hbox{s=c,t=0} \\
          \end{array}
  \right.
$$
\\
$$\leq \left\{
         \begin{array}{ll}
           \mathbb{E}\ [\ |h_{1 2\ 3 \ldots m}|^{\frac{5}{3}}\ \textbf{1}_{(|h|>n^{\frac{3s}{5}})} \ ]+\ \mathbb{E}\ [\ |h_{1 2 \ m+1\ldots 2m-2}|^{\frac{5}{3}}\ \textbf{1}_{(|h|>n^{\frac{3t}{5}})} \ ], & \hbox{s,t $>$ 0,\ s+t=c;} \\
          \mathbb{E}\ [\ |h_{1 2\ldots m}|^{\frac{5}{3}}\ \textbf{1}_{(|h|>n^{\frac{3c}{5}})} \ ]+\ \mathbb{P}\ (\ |h^{(m)}_{1 m+1 \ldots 2m-2}|>\ \log(n) \ ) , & \hbox{s=c,t=0}
         \end{array}
       \right.
$$
\\
$  \longrightarrow 0.$
\\ \\
Applying Markov inequality followed by an application Lemma 1 once again yields
\\
$$\mathbb{P}\ (\ |\ [n]^{-2m+2} \sum_{1\leq i_{1}\neq \ldots \neq i_{2m-2}  \leq }  \ V^{*{T^{\prime}}}(i_{k_{1}},\ldots,i_{k_{c}})\ |>\epsilon\ n\ (n-2m+2)\  )\qquad \qquad \qquad \qquad \qquad \qquad$$
\\
$$\leq K \epsilon^{-2}\ (n-2m+2)^{-2}  \ [n-(2m-2-c)]^{-c}\ n^{c+2}\ n^{-c-2}\ \mathbb{E}\ [\ h^{(s)}_{1 2  3\ldots m}\ h^{(t)}_{1 2\ m+1 \ldots 2m-2 }\ ]^{2}\qquad \qquad \qquad \ \ $$
\\
$$\leq \left\{
         \begin{array}{ll}
           K \epsilon^{-2}\ (n-2m+2)^{-2} \ [n-(2m-2-c)]^{-c}\ n^{c+2}\ n^{-c-2} \ n^{\frac{7c}{10}}\ \mathbb{E}|h_{1 2 \ldots m }|^{\frac{5}{3}}, & \hbox{s,t $>$ 0, s+t=c;} \\
           K \epsilon^{-2}\  (n-2m+2)^{-2} \ [n-(2m-2-c)]^{-c}\ n^{c+2}\ n^{-c-2}\ n^{\frac{7c}{10}}\ \log^{\frac{7}{6}}(n)\ \mathbb{E}|h_{1 2 \ldots m }|^{\frac{5}{3}} , & \hbox{s=c,t=0}
         \end{array}
       \right.
$$
$\longrightarrow 0,$ as $n\rightarrow\infty.$
\\
\\
This completes the proof of Proposition 4.2.
\\
\\
As the last step of the proof of Proposition 4, in the next result we  deal with terms of the form of sums of i.i.d. random variables (cf. Remark 7).
\\
\\
\textbf{\emph{Proposition 4.3}}.  \emph{If} $\mathbb{E}|h_{1\ldots m}|^{\frac{5}{3}}<\infty$,  \emph{then, as} $n\rightarrow \infty$
\\ \\
(a) $[n]^{-2m+1}\sum_{1\leq \neq i_{1}\neq \ldots \neq i_{2m-2} \leq n}V^{*}(i_{k_{1}})=o_{P}(1),\ k_{1}\in\{3,\ldots,2m-2\},$
\\
$$(b)\ \frac{1}{(n-2m+2) (n-2m+3)}\ \sum_{i\in\{1,\ldots,n\}/\{1,3,\ldots,2m-2\}}^{n} \mathbb{E}(h^{**}_{1 i 3 \ldots 2m-2}-\mathbb{E}(h^{**}_{1 i 3 \ldots 2m-2})\ |X_{i})=o_{P}(1), $$
\\
$$(c)\frac{1}{(n-2m+2) (n-2m+3)}\ \sum_{i\in\{1,\ldots,n\}/\{2,\ldots,2m-2\}}^{n} \mathbb{E}(h^{**}_{i 2 \ldots 2m-2}-\mathbb{E}(h^{**}_{i 2 \ldots 2m-2})\ |X_{i}) $$
$ \ \ \ +\ \displaystyle{\frac{1}{n-2m+2}}\ \mathbb{E}(h^{**}_{12\ldots 2m-2})\ =o_{P}(1).$
\\
\\
\textbf{\emph{Proof of Proposition 4.3}}
\\
\\
First we give the proof of part (a). Due to similarities, we shall state the proof only for the case that $k_{1}=3$.
\par
Define
\\
$$V^{*^{T}}(i_{3})=\ \mathbb{E}(h^{**^{T}}_{i_{1} i_{2}\ldots i_{2m-2}}- \mathbb{E}(h^{**^{T}}_{i_{1} i_{2}\ldots i_{2m-2}})|\ X_{i_{3}} ),$$
$$V^{*^{T^{\prime}}}(i_{3})=\ \mathbb{E}(h^{**^{T^{\prime}}}_{i_{1} i_{2}\ldots i_{2m-2}}- \mathbb{E}(h^{**^{T^{\prime}}}_{i_{1} i_{2}\ldots i_{2m-2}})|\ X_{i_{3}} ),$$
where $h^{**^{T}}_{i_{1} i_{2} \ldots i_{2m-2}}=\ h^{(1)}_{i_{1} i_{2} i_{3}\ldots i_{m}}\ h^{(m)}_{i_{1} i_{2} i_{m+1}\ldots i_{2m-2}}$ and $h^{**^{T^{\prime}}}_{i_{1} i_{2} \ldots i_{2m-2}}=\ h^{(1)}_{i_{1} i_{2} i_{3}\ldots i_{m}}\ h^{(0)}_{i_{1} i_{2} i_{m+1}\ldots i_{2m-2}}.$
\\
Again observe that as $n\rightarrow \infty $
\\
$$\mathbb{P}(\sum_{1\leq i_{1}\neq \ldots \neq i_{2m-2} \leq n}V^{*}(i_{3})\neq \sum_{1\leq i_{1}\neq \ldots \neq i_{2m-2} \leq n}V^{*^{T^{\prime}}}(i_{3})\ )\qquad \qquad \qquad \qquad \qquad \qquad $$
$$\leq\ \mathbb{P}(\sum_{1\leq i_{1}\neq \ldots \neq i_{2m-2} \leq n}V^{*}(i_{3})\neq \sum_{1\leq i_{1}\neq \ldots \neq i_{2m-2} \leq n}V^{*^{T}}(i_{3})\ )\qquad \qquad \qquad \qquad \qquad $$
$$+\ \mathbb{P}(\sum_{1\leq i_{1}\neq \ldots \neq i_{2m-2} \leq n} V^{*^{T}}(i_{3})\neq \sum_{1\leq i_{1}\neq \ldots \neq i_{2m-2} \leq n}V^{*^{T^{\prime}}}(i_{3})\ )  \qquad \qquad  \qquad \qquad $$
$$\leq\ n\ \mathbb{P}(|h^{(m)}_{1 2 3\ldots m}|>n^{3/5})+\ \mathbb{P}(|h^{(m)}_{1 2 \ m+1\ldots 2m-2}|>\log(n)) \qquad \qquad \qquad \qquad \qquad \qquad  $$
$\qquad \longrightarrow 0.$
\\
\\
Applying Markov inequality we arrive at
$$\mathbb{P}(\ |\frac{1}{n-2m+3} \sum_{1\leq i_{1}\neq \ldots \neq i_{2m-2} \leq n}V^{*^{T^{\prime}}}(i_{3})\ |>\ \epsilon\ (n-2m+2)  \ )\qquad \qquad \qquad \qquad $$
$$\leq \ K \ \epsilon^{-2}\ (n-2m+2)^{-2} \ (n-2m+3)^{-1}\ \mathbb{E}(h^{(1)}_{1 2 3\ldots m}h^{(0)}_{1 2 \ m+1\ldots 2m-2})^{2} \qquad \qquad $$

$\longrightarrow 0, \ \textrm{as}\  n\rightarrow\infty.$
\\
This completes the proof of part (a).
\par
Next to prove part (b) define
\\ \\
$h^{**^{T}}=\ h^{(1)}_{123\ldots m }h^{(m)}_{12\ m+1 \ldots 2m-2},$
\\
$h^{**^{T^{\prime}}}=\ h^{(1)}_{123\ldots m }h^{(1)}_{12\ m+1 \ldots 2m-2},$
\\ \\
and observe that as $n\rightarrow \infty$ we  have
$$\mathbb{P}\textbf{(}\sum_{i\in\{1,\ldots,n\}/\{1,3,\ldots,2m-2\}}^{n} \mathbb{E}(h^{**}_{1 i 3 \ldots 2m-2}-\mathbb{E}(h^{**}_{1 i 3 \ldots 2m-2})\ |X_{i})\qquad \qquad \qquad \qquad \qquad $$
$$\qquad \neq \sum_{i\in\{1,\ldots,n\}/\{1,3,\ldots,2m-2\}}^{n} \mathbb{E}(h^{**^{T^{\prime}}}_{1 i 3 \ldots 2m-2}-\mathbb{E}(h^{**^{T^{\prime}}}_{1 i 3 \ldots 2m-2})\ |X_{i})\ \textbf{)}\qquad \qquad \qquad \qquad  \qquad$$
$$\large{\leq}\ \mathbb{P}\textbf{(}\sum_{i\in\{1,\ldots,n\}/\{1,3,\ldots,2m-2\}}^{n} \mathbb{E}(h^{**}_{1 i 3 \ldots 2m-2}-\mathbb{E}(h^{**}_{1 i 3 \ldots 2m-2})\ |X_{i})\qquad \qquad \qquad \qquad $$
$$ \neq \sum_{i\in\{1,\ldots,n\}/\{1,3,\ldots,2m-2\}}^{n} \mathbb{E}(h^{**^{T}}_{1 i 3 \ldots 2m-2}-\mathbb{E}(h^{**^{T}}_{1 i 3 \ldots 2m-2})\ |X_{i})\ \textbf{)}\qquad \qquad \ \ \ $$
$$+\ \mathbb{P}\textbf{(}\sum_{i\in\{1,\ldots,n\}/\{1,3,\ldots,2m-2\}}^{n} \mathbb{E}(h^{**^{T}}_{1 i 3 \ldots 2m-2}-\mathbb{E}(h^{**^{T}}_{1 i 3 \ldots 2m-2})\ |X_{i})\qquad \qquad \qquad \qquad \qquad $$
$$\qquad \neq \sum_{i\in\{1,\ldots,n\}/\{1,3,\ldots,2m-2\}}^{n} \mathbb{E}(h^{**^{T^{\prime}}}_{1 i 3 \ldots 2m-2}-\mathbb{E}(h^{**^{T^{\prime}}}_{1 i 3 \ldots 2m-2})\ |X_{i})\ \textbf{)}  \qquad \qquad \qquad \qquad \ $$
$$\large{\leq}\ n\ \mathbb{P}(|h^{(m)}_{1 2 3\ldots m}|>n^{3/5})+ \  n\ \mathbb{P}(|h^{(m)}_{1 2 \ m+1\ldots 2m-2}|>n^{3/5})\qquad \qquad \qquad \qquad$$
$\longrightarrow 0.$
\\ \\
Hence another application of Markov inequality yields,
$$\mathbb{P}\textbf{(}\ |\frac{1}{n-2m+3} \sum_{i\in\{1,\ldots,n\}/\{1,3,\ldots,2m-2\}}^{n} \mathbb{E}(h^{**^{T^{\prime}}}_{1 i 3 \ldots 2m-2}-\mathbb{E}(h^{**^{T^{\prime}}}_{1 i 3 \ldots 2m-2})\ |X_{i})\  | > \epsilon \ (n-2m+2)   \ \textbf{)}$$
$$\leq \ K \ \epsilon^{-2}\ (n-2m+2)^{-2} \ (n-2m+3)^{-1}\ \mathbb{E}(h^{(1)}_{1 2 3\ldots m}h^{(1)}_{1 2 \ m+1\ldots 2m-2})^{2} \qquad \qquad $$
$\longrightarrow 0, \ \textrm{as} \ n\rightarrow\infty.$
\\
Now the proof of part (b) is complete.
\\
\par
To prove part (c) we only need to observe that
$$\frac{1}{(n-2m+2) (n-2m+3)}\ \sum_{i\in\{1,\ldots,n\}/\{2,\ldots,2m-2\}}^{n} \mathbb{E}(h^{**}_{i 2 \ldots 2m-2}-\mathbb{E}(h^{**}_{i 2 \ldots 2m-2})\ |X_{i}) $$
$$ \ + \displaystyle{\frac{1}{n-2m+2}}\ \mathbb{E}(h^{**}_{12\ldots 2m-2})= \ \frac{1}{(n-2m+2) (n-2m+3)}\ \sum_{i\in\{1,\ldots,n\}/\{2,\ldots,2m-2\}}^{n} \mathbb{E}(h^{**}_{i 2 \ldots 2m-2}\ |X_{i}).$$
The rest of the proof is similar to that of part (b), hence the details are omitted.
Now the proof of Proposition 4.3 and that of (I) is complete.
\\

\textbf{Acknowledgments.} The author wishes to thank Mikl\'{o}s Cs\"{o}rg\H{o}, Barbara Szyszkowicz and Qiying Wang for calling his attention to a preliminary version of their paper \cite{cs8b} that inspired the truncation arguments of the present exposition. This work constitutes a part of the author's Ph.D. thesis in preparation, written under the supervision and guidance of Mikl\'{o}s Cs\"{o}rg\H{o} and Majid Mojirsheibani. My special thanks to them for also reading preliminary versions of this article, and for their instructive comments and suggestions that have much improved the construction and presentation of the results.


\begin{thebibliography}{9999}
\bibitem{Arv} Arvesen, J. N. (1969). Jackknifing U-statistics. \emph{Ann. Math. Stat.} \textbf{40 }, 2076-2100
\bibitem{cs3} Cs\"{o}rg\H{o}, M., Szyszkowicz, B. and Wang, Q. (2003). Donsker's theorem for self-normalized parial sums processes. \emph{The Annals of Probability} \textbf{31}, 1228-1240.
\bibitem{cs4} Cs\"{o}rg\H{o}, M., Szyszkowicz, B. and Wang, Q. (2004). On Weighted
Approximations and Strong Limit Theorems for Self-normalized
Partial Sums Processes. In\emph{ Asymptotic methods in Stochastics},
489-521, Fields Inst. Commun.44, Amer. Math. Soc., Providence RI.
\bibitem{cs8a}Cs\"{o}rg\H{o}, M., Szyszkowicz, B. and Wang, Q. (2008). On weighted approximations in $D[0,1]$ with application to self-normalized partial sum processes.
\emph{Acta Mathematica Hungarica} \textbf{121 (4)}, 307-332.
\bibitem{cs8b} Cs\"{o}rg\H{o}, M., Szyszkowicz, B. and Wang, Q. (2008). Asymptotics of studentized U-type processes for changepoint problems. \emph{Acta Mathematica Hungarica} \textbf{121 (4)}, 333-357.
\bibitem{gin1} Gin\'{e}, E., and Zinn, J. (1992). Marcinkiewicz type laws of large numbers and convergence
of moments for U-statistics. In \emph{Probability in Banach Spaces} (R. Dudley, M. Hahn
and J. Kuelbs, eds) \textbf{8} 273-291, Birkhauser, Boston.
\bibitem{gin2} Gin\'{e}, E. , G\"{o}tze, F. and Mason D. M. (1997). When is the student t-statistic asymptotically Normal? \emph{The Annals of Probability} \textbf{25}, 1514-1531.
\bibitem{gut} Gut, A. (2005). \emph{Probability: A Graduate Course.} Springer.
\bibitem{hof} Hoeffding, W. (1948). A class of statistics with asymptotically normal distribution. \emph{Ann. Math. Statist}. \textbf{19}, 293-325.
\bibitem{mil} Miller, R. G. Jr. and Sen, P. K. (1972). Weak convergence of U-statistics and Von Mises' differentiable statistical functions. \emph{Ann. Math. Statist.} \textbf{43}, 31-41.
\bibitem{nasa} Nasari, M. M. (2009). On weak approximations of $U$-statistics. \emph{Statistics and Probability letters,} \textbf{79}, 1528-1535.
\bibitem{serf}Serfling, R. J. (1980). \emph{Approximation Theorems of Mathematical
Statistics}. Wiley, New York.
\end{thebibliography}
\end{document}